\documentclass[12pt,onecolumn]{IEEEtran}
\usepackage[margin=2.5cm]{geometry}
\usepackage[cmex10]{amsmath}

\usepackage{mathtools}

\usepackage[tight,footnotesize]{subfigure}
\usepackage[caption=false,font=footnotesize]{subfig}
\usepackage{bm}
\usepackage{amsmath}
\usepackage{amssymb}
\usepackage{epsfig}
\usepackage{wrapfig}
\usepackage{array}
\usepackage{engord}
\usepackage{amsthm}
\usepackage{graphicx}
\usepackage{epstopdf}
\usepackage{siunitx}

\newcommand{\hypone}[3]{\ensuremath{H_{#1, #2}(#3)}}
\newcommand{\hypk}[4]{\ensuremath{H_{#1, #2}(#3, #4)}}

\newcommand{\hypinit}[2]{\ensuremath{\hypone{1}{#1-1}{#2}}}

\newcommand{\riskone}[3]{\ensuremath{R_{#1, #2}(#3)}}
\newcommand{\riskk}[4]{\ensuremath{R_{#1, #2}(#3, #4)}}
\newcommand{\dval}[3]{\ensuremath{d_{#1}(#2,#3)}}

\newcommand{\lkk}[8]{\ensuremath{L(\hypk{#1}{#2}{#3}{#4},\hypk{#5}{#6}{#7}{#8})}}
\newcommand{\lkone}[7]{\ensuremath{L(\hypk{#1}{#2}{#3}{#4},\hypone{#5}{#6}{#7})}}
\newcommand{\lonek}[7]{\ensuremath{L(\hypone{#1}{#2}{#3},\hypk{#4}{#5}{#6}{#7})}}
\newcommand{\loneone}[6]{\ensuremath{L(\hypone{#1}{#2}{#3},\hypone{#4}{#5}{#6})}}
\newcommand{\setsquare}{\ensuremath{ \mathcal{S}^{\overline{2}} }}

\graphicspath{{./graphics/}{.}}
\usepackage{color}
\usepackage{url}
\usepackage{cite}
\usepackage{verbatim}

\newtheorem{theorem}{Theorem}

\newtheorem{remark}{Remark}
\newtheorem{lemma}{Lemma}
 \IEEEoverridecommandlockouts
\begin{document}
\title{Sequential Detection of an Abrupt Change in a Random Sequence with Unknown Initial State}
\author{\IEEEauthorblockN{James Falt and Steven~D.~Blostein, \IEEEmembership{Senior Member,~IEEE}\thanks{The authors are with Dept.\ of Electrical and Computer Eng., Queen's University, Kingston, ON, Canada. (contact email: steven.blostein@queensu.ca).  Early versions of this  paper were presented in part in \cite{Falt2014} and \cite{Falt2016}. This research was supported by NSERC Discovery Grant 05061-2014. }}}
\maketitle
\thispagestyle{empty}
\begin{abstract} 
The problem of sequentially detecting an abrupt change in a sequence of independent and identically distributed (IID) random variables is addressed. 
Whereas previous approaches assume a known probability density function (PDF) at the start of the sequence, the problem addressed is the detection of a single change in distribution among a finite number of known 'equal-energy' PDFs, but  where the initial and final distributions are not known a priori. 
A  Bayesian multiple hypothesis approach is proposed where (i) unlike previous threshold policies,  the minimum cost hypothesis is tracked through time, (ii) under an exponential delay-cost function that satisfies an upper bound determined by the distances between hypotheses, the probability of detecting a change from an incorrect initial distribution asymptotically vanishes with time, (iii)  computation is recursive and constant per unit time, and (iv) the unknown initial state gives rise to unavoidable incorrect detections that be made to vanish with a constant test threshold with negligibly small effect on correct detection delay for change times beyond a lower bound. Simulations illustrate the analysis and reveal that average delay approaches that of the optimal CUSUM test after an initial transient period determined by an incorrect detection probability constraint.
\end{abstract}

\section{Introduction}

Let $ X_{1,n} = \{ X_i | i = 1,2,\dots,n \} $ be $n$ independent random variables observed sequentially. Each of these random variables are known to follow one of $D$ possible known distributions, which are described by probability density functions (PDF) $f_j$, for $j \in \{ 0, 1, \dots, D-1\}$, each with finite variance. While observing $X_{1,n}$ sequentially, the objective is to determine, as quickly as possible, whether a single change in distribution has occurred at some unknown discrete time $m \in \mathcal{Z+}$ \emph{without prior knowledge of the initial or final distributions}. Generalizations of the standard ($D=2$) quickest detection problem formulation have been considered in \cite{Blostein1991,Moustakides1998} but for only limited cases. Specifically, \cite{Blostein1991} addresses the case where the distribution of the sequence is time-varying after the change time $m$, and \cite{Moustakides1998} considers a random process which is not necessarily IID before and after the change. In both cases, however, the initial PDF  is assumed to be known.

Some of change detection's traditional applications can be found in \cite{Basseville93}. An example of more recent interest  that features an unknown initial state with $D=2$ arises in cognitive radio,  where $f_0$ represents the PDF of noise only, representing the absence of a coherent pilot signal denoting usage of a communication channel,  while  $f_1$ represents the pilot's presence in the same background noise \cite{Liang11,IEEE802.22Introduction}.  
In contrast to previous formulations, the {\em initial state} corresponding the PDF at the start of observation $X_1$ is unknown. In addition, conditioned on either $f_0$ or $f_1$, observations are assumed to be independent.

In the case where the initial state is known, solutions to the problem of detection of a single change in distribution in an independent sequence of random variables can be found in\cite{Page1954,Lorden1971,Shiryaev1978,Moustakides1986}, as well as the substantial number of references that can be found in  \cite{Poor2008}. More recently, the problem of detecting multiple changes in sequences, termed {\em ON-OFF processes} has been addressed in \cite{Zhao2010}.
Nevertheless, as in previous approaches, the initial state in \cite{Zhao2010} is assumed known.  It is also worth noting that quickest detection problems yield threshold policies, and the problem of threshold determination to guarantee an exact probability of false alarm has been generally intractable except for a few special cases, e.g., involving a priori Bernoulli PDFs and geometric change time distribution \cite{Shiryaev1978,Shiryaev1963,Poor2008}.

In the literature there is also what is termed a {\em 2-sided test,} i.e., a sequence having two possible final states, as first proposed in \cite{Barnard1959}, analyzed asymptotically and for Normal distributions in\cite{Dragalin1997} and applied to certain problems \cite{Barnett2001}, as well as the detection of changes in drift of Brownian motion  \cite{Hadjiliadis2005,Hadjiliadis2006,Hadjiliadis2008}. In the latter problem, asymptotic optimality is proven as the frequency of false alarms tends to infinity as well as in an asymptotic minimax sense. In all of the above, although the change or drift itself may be 2-sided, the initial state of the process is assumed completely known. A decentralized version of the 2-sided, known initial-state  problem is considered in \cite{Zarrin09}. However, it has been observed that quickest detection performance is drastically affected if the assumed initial state is incorrect \cite{Liang11}. 

%Change detectors, such as CUSUM \cite{Page1954} and Shiryaev  \cite{Shiryaev1978}, assume knowledge of initial state. Classical detection is concerned with two types of possible errors: false alarm and missed detection. A false alarm is defined as the event where a stopping rule is invoked before a change has occurred. A missed detection is defined as the event where a stopping rule is never invoked after a change occurs. 
In previous approaches to one-sided change detection where the initial and final states are known, different performance metrics are used  depending on the assumptions about prior information. Bayesian formulations assume that the change time is random with known prior distribution and  minimize average detection delay subject to a constraint on the probability of false alarm \cite{Shiryaev1963,Shiryaev1978}. Alternatively,  an unknown non-random change time is assumed and formulated in a minimax sense in \cite{Lorden1971}, where it has found that  the CUSUM procedure by Page \cite{Page1954} minimizes  worst-case detection delay  subject to an upper bound on the false alarm rate (FAR), or equivalently, a lower bound on the inverse of FAR, known as the average run length (ARL) to false alarm. 
%The metric characterizing the detection delay in Lorden's criterion is the worst-case detection delay over all possible change points, whereas in the Shiryaev-Roberts-Pollak criterion it is the worst-case average detection delay.
%Under this assumption, it is not possible to calculate the probability of false alarm, since this probability depends on the change time itself. 
%Tests which are formulated under this assumption, such as those of Page \cite{Page1954} and Shiryaev \cite{Shiryaev1978}, use the average detection delay and the average run length (ARL) to false alarms to characterize test performance since the probability of false alarm depends on the change time. 
More recently, the CUSUM procedure  has been shown to be optimal in terms of minimizing average delay itself rather than an upper bound 
\cite{Moustakides1986, Pollak85, Ritov90}.

In contrast to the above, the proposed approach to change detection with unknown initial state follows a Bayesian formulation that minimizes expected delay for a specific proposed cost structure. 
%In the proposed change detector's formulation it was assumed that the change time is unknown, and thus average detection delay for a deterministic change time and the FAR will be used to characterize the test's performance.
%In previous approaches to change detection where the starting state is assumed to be known, such as those of CUSUM \cite{Page1954} and that of Shiryaev  \cite{Shiryaev1978}, the probability of a missed detection can easily be made to vanish. Rather than missed detection, average detection delay has therefore been used to assess change detector performance. 
Since initial and final PDFs are unknown, we define {\em incorrect detection} of events whenever the initial and/or final states are incorrectly declared.
%to indicate the event where a change is declared and the initial state chosen is incorrectly. Additionally, for $D>2$, we must also consider the detection error where the initial state is identified correctly but the final state is chosen incorrectly following the change in distribution.
%We note that for the events $H_1$ and $\overline{H_1}$ incorrect detection 
Since information is being gathered over time,
%as the number of samples observed from the initial state increases, 
we seek a procedure where the probability of incorrect detection of initial state decreases with time.
A sufficient condition is if the Bayes' risks associated with hypotheses with correct initial state converge, while risks corresponding to hypotheses with incorrect initial state become unbounded.
In other words, once enough samples are observed, we seek a test with performance that approaches that of a change detector with a known initial state and one of possibly multiple known final states.

This paper extends, generalizes and builds on earlier versions \cite{Falt2014} and \cite{Falt2016}. A Bayesian change detector is first formulated in \cite{Falt2014} under unknown initial state with an exponential delay-cost model, but restricted to two states, i.e, $D=2$, with only asymptotic behaviour considered, leading to a preliminary version of Theorem~1 in the sequel. The model was improved in \cite{Falt2016} to include initial state uncertainty cost, thereby lowering the probability of incorrect initial state of the sequence. Here, we generalize \cite{Falt2016}, to allow for $D \geq 2$ distributions with equal-energy finite-variance PDFs.  In addition, the performance analysis and results are extended, strengthened and more fully described.

While uniform costs and equal priors are sufficient for determining the most likely hypothesis, the causal nature of the problem requires detection delays to be penalized. In Section~\ref{sec-2sided}, the change detection problem for unknown initial state is therefore formulated under an exponential delay-cost model. Exponential delay cost has been proposed elsewhere \cite{Poor1998}. Section~\ref{sec-2sided} proposes a recursive formulation  suitable for problems involving causal (on-line) observation, that  implemented by tracking a fixed number of states, which is  a consequence of the chosen cost structure.
While full joint optimization of the test's multiple parameters is intractable, Section~\ref{sec-design} presents bounds and tradeoffs that govern suitable parameter choices.
%While design of the test threshold is also intractable here,  a parameterization of the cost function is proposed and studied in Section~\ref{sec-design}. 
Section~\ref{sec-results} presents Monte Carlo simulations with insights on the performance analysis of Section ~\ref{sec-design}, and draws performance comparisons to optimal one-sided quickest detection methods for known initial state\cite{Moustakides1986} \cite{Poor2008}. 
%Results are described in Section~\ref{sec-results} which reinforce the insights provided in Section~\ref{sec-design} as well as make comparisons to CUSUM, which has been shown to be optimal for one-sided quickest detection \cite{Moustakides1986} \cite{Poor2008}. 

\section{Change Detection under Unknown Initial State: Bayes Formulation} \label{sec-2sided}

At time $n$, without knowledge of the starting PDF (of $X_1$), a sequence of independent random variables, $X_{1,n} \equiv  X_1 , X_2, \cdots X_n$, is observed. The initial distribution of the sequence is one of $D$ possible distributions, which are given by the distinct {\em equal-energy} PDFs $f_i$, for $i \in \{ 0, 1, \dots, D-1\}$, where $\mathbb{E}_i [f_i(X)] = \mathbb{E}_j [f_j(X)]$ for $i,j \in \{ 0, 1, \dots, D-1\}$. The sequence may assume $D + (n-1)\frac{D!}{(D-2)!}$ possible joint distributions, of which 
%can be enumerated as follows: 
$D$ of the sequences, corresponding to {\em no change}, are described by the sequence $X_{1,n}$ where all $ X_i \sim f_j, 1 \leq i \leq n$, where 
$j \in \mathcal{S} = \{0,1, \dots, D-1\}$. The remaining $(n-1)\frac{D!}{(D-2)!}$ possible sequences, which correspond to each of the possible single distribution changes which occur after the first sample, are described by the sequence $X_{1,n}$ where $X_i  \sim f_j , 1 \leq i < m$, and the rest of the $X_i  \sim f_k, m \leq i \leq n$, where $(j,k) \in \setsquare = \{ (j,k) \; | \;  j,k \in \mathcal{S}  \mbox{ \small{and} }  j \neq k \}$ and $1 < m \leq n$ is the change time. At time $n$, there are $(n-1)\frac{D!}{(D-2)!}$ change hypotheses as there are $n-1$ possible change times after the first sample and $\frac{D!}{(D-2)!}$ ordered subsets of 2 elements in $\mathcal{S}$. Classically, for a fixed number of samples, the detection problem is one of testing a finite number of hypotheses \cite{Poor1994} given by $M(n) \equiv D + (n-1)\frac{D!}{(D-2)!}$.
%$ j,k  \; \in  \; \{ 0,1, \dots, M-1 \}$, $i \neq j$, and $m$ is the {\em change time}. 
The notation adopted to represent these hypotheses is the following:
\begin{eqnarray} \hypone{1}{n}{j} &=& \{  X_i \sim f_j |  1 \leq i \leq n\}  \qquad j \in \mathcal{S} \label{h1}\\
			% \overline{H_1} &=& \{  X_i \sim f_1 |  1 \leq i \leq n\} \label{h1overline} \\
			\hypk{m}{n}{j}{k} &=& \left\{ \begin{array}{ll} X_i \sim f_j | 1 \leq i < m \\
								X_i \sim f_k | m \leq i \leq n
						\end{array} \right\} \qquad 1 < m \leq n, \; (j,k) \in \setsquare .\label{hj}
			% \overline{H_j} &=& \left\{ \begin{array}{ll} X_i \sim f_1 | 1 \leq i < j \\
			% 					X_i \sim f_0 | j \leq i \leq n
			% 			\end{array} \right\} \qquad 2 \leq j \leq n \label{hjoverline}.
\end{eqnarray}
In the above, when $\hypone{1}{n}{j}, \; j \in \mathcal{S}$ is selected, there is {\em no detection} of a change.  If $\hypk{m}{n}{j}{k},$ $ \; (j,k) \in \setsquare, \; 1 < m \leq n,$  is selected before a change occurs when $\hypone{1}{n}{j}$ is true,  a {\em false alarm} occurs. On the other hand,  if $\hypone{1}{n}{j}, \; j \in \mathcal{S}$ is selected and any change has occurred at time $1 < m < n$ there is {\em detection delay}.  If $\hypk{m}{n}{j}{k}$, for $1 < m \leq n$ is selected, and the true hypothesis is $\hypk{l}{n}{j}{s}$, for $1 < l \leq n$, $(j,k) \in \setsquare$, $(j,s) \in \setsquare$, and $k \neq s$, then an {\em incorrect detection of final state} occurs. If $\hypk{m}{n}{j}{k}$, for $1 < m \leq n$ is selected, and the true hypothesis does not have $f_j$ as its initial distribution, then an {\em incorrect detection from initial state} occurs, regardless of the accuracy of the selected final state $f_k$ or the selected change time $m$.  

Suppose for now that time $n$ is fixed. To formulate a decision rule to select among the $M(n)$ possible hypotheses, we adopt a Bayesian framework.
% based on selecting the hypothesis with minimum risk.
% according to Bayesian $M(n)$-ary hypothesis testing. 
Under equally likely prior probabilities of change time,  the Bayesian test minimizes the posterior risk associated with each of the $M(n)$ hypotheses. However, when an infinite-duration sequence is causally observed,  the number of hypotheses grows with $n$, which results in increasing computation and memory. To alleviate this, a time-recursive version of Bayes risk computation testing is required. 

The posterior probabilities for each hypothesis  given $X_{1,n} = x_{1,n}$ need to be tracked over time, and we denote the posterior probability that the hypothesis $\hypk{m}{n}{j}{k}$, $ (j,k) \in \setsquare$ is true given $X_{1,n} = x_{1,n}$  by $P( \hypk{m}{n}{j}{k} | X_{1,n} = x_{1,n} )$ for $ 1 < m \leq n $, and similarly the posterior probability that the hypothesis $\hypone{1}{n}{j} , \; j \in \mathcal{S}$ is true given $X_{1,n} = x_{1,n}$ is denoted by $P( \hypone{1}{n}{j} | X_{1,n} = x_{1,n} )$. Using Bayes' Rule,  
\begin{eqnarray}
P( \hypk{m}{n}{j}{k} | X_{1,n} = x_{1,n} ) &=& \frac{P( X_{1,n} = x_{1,n} | \hypk{m}{n}{j}{k} )   P(\hypk{m}{n}{j}{k})}{P(X_{1,n} = x_{1,n})} \label{bayes}
\end{eqnarray}
where $P( X_{1,n} = x_{1,n} | \hypk{m}{n}{j}{k} )$ is the likelihood of observing $X_{1,n} = x_{1,n}$ given that $\hypk{m}{n}{j}{k}$ is true, $P(\hypk{m}{n}{j}{k})$ is the prior probability of $\hypk{m}{n}{j}{k}$, and $P(X_{1,n} = x_{1,n})$ is the likelihood of observing $X_{1,n} = x_{1,n}$ for the $n$ samples received. An equivalent expression to (\ref{bayes}) can be written for the posterior probability of the hypothesis $\hypone{1}{n}{j} , \; j \in \mathcal{S}$  given $X_{1,n} = x_{1,n}$. From the assumed independence conditioned on a certain hypothesis, the likelihood functions  in (\ref{bayes}) are 
\begin{align}
&P( X_{1,n} = x_{1,n} | \hypone{1}{n}{j} ) = \prod_{i=1}^{n} f_j(x_i) && j \in \mathcal{S}\label{posterior_one} \\
&P( X_{1,n} = x_{1,n} | \hypk{m}{n}{j}{k} ) = \prod_{i=1}^{m-1} f_j(x_i) \prod_{i=m}^{n} f_k(x_i) && (j,k) \in \setsquare, \; 1 < m \leq n.\label{posterior_k}
\end{align}
Define hypotheses $H1$ and $H2$, each of which assume the form of either (\ref{h1}) or (\ref{hj}).
The proposed Bayes formulation uses cost function $L(H1,H2)$, which is the cost of choosing $H1$ when $H2$ is true, 
%\begin{eqnarray*}
%\lkk{m_1}{n}{j_1}{k_1}{m_2}{n}{j_2}{k_2} &\equiv& \mbox{ cost of choosing $\hypk{m_1}{n}{j_1}{k_1}$ when $\hypk{m_2}{n}{j_2}{k_2}$ is true} \\% , $(j_1,k_1), (j_2,k_2) \in Z^2$ }\\
%\lonek{1}{n}{j_1}{m_2}{n}{j_2}{k_2} &\equiv& \mbox{ cost of choosing $\hypone{1}{n}{j_1}$ when $\hypk{m_2}{n}{j_2}{k_2}$ is true}\\%, $j_1 \in Z, \; (j_2,k_2) \in Z^2$ } \\
%\lkone{m_1}{n}{j_1}{k_1}{1}{n}{j_2} &\equiv& \mbox{ cost of choosing $\hypk{m_1}{n}{j_1}{k_1}$ when $\hypone{1}{n}{j_2}$ is true}\\%, $(j_1,k_1) \in Z^2, \; j_2 \in Z$ }\\
%\loneone{1}{n}{j_1}{1}{n}{j_2} &\equiv& \mbox{ cost of choosing $\hypone{1}{n}{j_1}$ when $\hypone{1}{n}{j_2}$ is true,}%, $j_1, j_2 \in Z$}
%\end{eqnarray*}
%where $1 < m_1 \leq n$, $1 < m_2 \leq n$, and $(j_1,k_1),(j_2,k_2) \in \setsquare$. 
and the priors 
\begin{align*}
\pi_m{(j,k)} &\equiv \mbox{ prior probability that $\hypk{m}{n}{j}{k}$ is true, $(j,k) \in \setsquare,\; m > 1$}\\
\pi_1(j) &\equiv \mbox{ prior probability that $\hypone{1}{n}{j}$ is true, $j \in \mathcal{S}.$}
\end{align*}
While uniform costs $\loneone{1}{n}{j_1}{1}{n}{j_2} = 0$ for $j_1=j_2 \in \mathcal{S}$, $\lkk{m_1}{n}{j_1}{k_1}{m_2}{n}{j_2}{k_2}$$ = 0$ for $m_1=m_2$ and $(j_1,k_1) = (j_2,k_2) \in \setsquare$, and all costs otherwise equal 1 may be used to determine the maximum a posteriori hypothesis at each $n$,  alternative cost structures to reflect the problem's  time-sequential nature are more appropriate. Exponential cost has been explored in \cite{Poor1998} for change detection problems where the initial state is known, and allow for performance trade offs between average detection delay and false alarm rate. 

For change detection under unknown initial state, undesired incorrect detections need to be controlled. Specifically, incorrect detection arising from {\em initial state uncertainty} of the sequence $\{ X_i | i = 1,2,\cdots,n \} $ occurs if  a change from $f_j$ to $f_k$ is declared while $f_k$ is the initial state, for $(j,k) \in \setsquare$. Hypothesis $\hypk{m}{n}{j}{k}$, for $1 < m \leq n$, has its first $m-1$ samples correspond to the initial state. A non-sequential Bayesian test of all $M(n)$-ary possible hypotheses does not associate a cost with uncertainty in initial state, as the likelihood of a certain hypothesis is a function of all $n$ observations.
%since the posterior probabilities for each hypothesis are calculated using all of the $n$ observed samples, which for 
% $m<n$ will 
%not all are distributed according to the initial state. 
The notion of Bayes' risk for initial state uncertainty for each of the possible $(n-1)\frac{D!}{(D-2)!}$ change hypotheses is therefore needed in the formulation. Let $\hypinit{m}{j}$, $j \in \mathcal{S}$, denote the hypothesis with change time $1 < m \leq n$, initial state $f_j$, and prior probability $\phi(j)$. Using Bayes' rule, the posterior probabilities are 
\begin{eqnarray} \label{initposterior0}
P( \hypinit{m}{j} | X_{1,m-1} = x_{1,m-1} ) = \frac{P(  X_{1,m-1} = x_{1,m-1} | \hypinit{m}{j} ) P(\hypinit{m}{j}) }{P(  X_{1,m-1} = x_{1,m-1})}, \;j \in \mathcal{S}.
\end{eqnarray}
%\begin{eqnarray} \label{initposterior1}
%P( H_k^1 | X_1^{k-1} = x_1^{k-1} ) = \frac{P(  X_1^{k-1} = x_1^{k-1} | H_k^1 ) P(H_k^1)}{P(  X_1^{k-1} = x_1^{k-1})}.
%\end{eqnarray}
The likelihood term in (\ref{initposterior0}) factors using independence  as
\begin{eqnarray}
P(  X_{1,m-1} = x_{1,m-1} | \hypinit{m}{j}) &=& \prod_{i=1}^{m-1} f_j (x_i), \; j \in \mathcal{S}.
% P(  X_1^{k-1} = x_1^{k-1} | H_k^1 ) &=& \prod_{i=1}^{k-1} f_1 (x_i).
\end{eqnarray}
Define wrong initial state cost $I(j,k) \equiv \mbox{ cost of choosing} \; \hypinit{m}{j}  \mid  \hypinit{m}{k}  \; \mbox{true}, j,k \in  \mathcal{S}$.

The posterior risk in choosing hypothesis $\hypone{1}{n}{j}, \; j \in \mathcal{S}$, which corresponds to {\em no change}, given $X_{1,n} = x_{1,n}$, is 
\begin{align}
\riskone{1}{n}{j} = &  \sum_{\{ r \in \mathcal{S} \}} \Bigg{(} \loneone{1}{n}{j}{1}{n}{r} P( \hypone{1}{n}{r} | X_{1,n} = x_{1,n}) \nonumber \\
& \hspace{5pt}+ \sum_{\{ s \in \mathcal{S}| s \neq r \} }  \sum_{i=2}^{n} \lonek{1}{n}{j}{i}{n}{r}{s} P ( \hypk{i}{n}{r}{s} | X_{1,n} = x_{1,n} )  \Bigg{)}.
\label{R1n}
\end{align}

The posterior risk in choosing hypothesis $\hypk{m}{n}{j}{k}, \; 1 < m  \leq n, \; (j,k) \in \setsquare$, 
%given $X_{1,n} = x_{1,n}$, 
%can be expressed as 
\begin{align}
\riskk{m}{n}{j}{k} = & \sum_{ \{r \in \mathcal{S} \} } \Bigg{[}  \lkone{m}{n}{j}{k}{1}{n}{r} P( \hypone{1}{n}{r} | X_{1,n} = x_{1,n}) \nonumber \\
& \hspace{35pt}+ \sum_{ \{ s \in \mathcal{S}| s \neq r \} }  \sum_{i=2}^{n}  \lkk{m}{n}{j}{k}{i}{n}{r}{s} P ( \hypk{i}{n}{r}{s} | X_{1,n} = x_{1,n}   \Bigg{]} \nonumber \\
& + \sum_{\{ r \in \mathcal{S} \} }  I(j,r) P( \hypinit{m}{r} | X_{1,m-1} = x_{1,m-1}) .
\label{Rkn}
\end{align}

%The conditional risk in choosing hypothesis $H_k, \; 1 \leq k \leq n$,  is given $X_1^n = x_1^n$,  
%\begin{equation}
%R_k^n = \sum_{j=1}^n ( \pi_j   L_{kj}   P( H_j | X_1^n = x_1^n ) + \overline{\pi_j}   L_{k\overline{j}}   P( \overline{H_j} | X_1^n = x_1^n ) ), \label{Rkn}
%\end{equation}
%and the conditional risk in choosing $\overline{H_k}$ given $X_1^n = x_1^n$, is given by
%\begin{equation}
%R_{\overline{k}}^n = \sum_{j=1}^n ( \pi_j   L_{\overline{k}j}   P( H_j | X_1^n = x_1^n ) + \overline{\pi_j}   L_{\overline{k}\overline{j}}   P( \overline{H_j} | X_1^n = x_1^n )) . \label{koverlinenrisk}
%\end{equation}

The proposed stopping rule, $\psi (n)$,  is to stop at the first $n$ such that
\begin{align}
\psi (n):   \min \left\{  \riskk{m}{n}{j}{k} \; (j,k) \in \setsquare, \; 1 < m \leq n \right\} < \min   \left\{    \riskone{1}{n}{i},  i \in \mathcal{S} \right\}. \label{minrisk}
\end{align}

From here on equally likely a priori hypotheses, i.e., $\pi \equiv \pi_{1}{(j)} = \pi_{m}{(j,k)} = 1/M(n), \; (j,k) \in \setsquare$, and equally likely initial states, i.e.,  $\phi \equiv  \phi(j) = 1/D, \; j \in \mathcal{S}$ are assumed.
To penalize detection delay, false alarm, and incorrect detection, the following cost structure is adopted:
%In the cost $L(H1,H2)$, we define the hypothesis $H1$ as the \emph{subject} hypothesis, and $H2$ as the \emph{alternative} hypothesis. %If $H1$ and $H2$ are both change hypotheses, the hypotheses are said to be \emph{non-conflicting} if they both have the same initial and final states. If one or both of the hypotheses are no-change hypotheses, the hypotheses $H1$ and $H2$ are non-conflicting if they share the same initial state. Otherwise, the hypotheses $H1$ and $H2$ are said to be \emph{conflicting}. 
%The cost structure consists of four distinct types of costs $L(H1,H2)$:
\begin{enumerate}
\item \textbf{Zero Cost:} If $H1 = H2$, $L(H1,H2)=0$.
\item \textbf{Fixed Costs of False Alarm and Incorrect Final State:} If $H1$ and $H2$ share initial and final states and $H1$ corresponds to a change occurring earlier than in $H2$, $L(H1,H2)$ has a fixed cost of false alarm, $b$. Additionally, if $H1$ and $H2$ share initial states but not final states, $L(H1,H2)$ has a fixed cost of incorrect final state, which is also $b$.
\item \textbf{Exponential Cost of Delay:} If $H1$ and $H2$ share initial and final states and $H1$'s  change time is greater than that of $H2$, $L(H1,H2)$ has the exponential cost of delay, $a^i$, where $i$ is the delay by which $H1$ lags $H2$. 
\item \textbf{Exponential Cost of Incorrect Detection from Initial State:} If $H1$ and $H2$ have different initial states, $L(H1,H2)$ has the exponential cost of incorrect detection from initial state, $c^i$, where $i$ is  the number of samples in which $H1$ and $H2$ are distributed differently.
\end{enumerate}

%
%
%Cost Structure Explicitly Here%

The cost structure described above is explicitly given by:
\begin{align}
\lkk{m_1}{n}{j_1}{k_1}{m_2}{n}{j_2}{k_2} &\equiv \left\{ \begin{array}{ll} 
						\left. \begin{array}{ll}
							b & \mbox{if $j_1 = j_2$ and $k_1 = k_2$ } \\
							b & \mbox{if $j_1 = j_2$ and $k_1 \neq k_2$ } \\
							c^{m_2-1} & \mbox{if $j_1 \neq j_2$ and $k_1 = k_2$ }  \\
							c^{n-m_2+1}c^{m_1-1} & \mbox{if $j_1 \neq j_2$, $k_1 \neq k_2$, and $k_1 = j_2$ } \\
							c^n & \mbox{if $j_1 \neq j_2$, $k_1 \neq k_2$, and $k_1 \neq j_2$}
						\end{array} \right\} m_2 > m_1 \\
						\left. \begin{array}{ll}
							a^{m_1-m_2} & \mbox{if $j_1 = j_2$ and $k_1 = k_2$ } \\
							b & \mbox{if $j_1 = j_2$ and $k_1 \neq k_2$ } \\
							c^{m_1-1} & \mbox{if $j_1 \neq j_2$ and $k_1 = k_2$ }  \\
							c^{n-m_1+1}c^{m_2-1} & \mbox{if $j_1 \neq j_2$, $k_1 \neq k_2$, and $k_2 = j_1$ } \\
							c^n & \mbox{if $j_1 \neq j_2$, $k_1 \neq k_2$, and $k_2 \neq j_1$}
						\end{array} \right\} m_1 > m_2 \\
						\left. \begin{array}{ll}
							b \;\;\;\;\;\;\;\;\;\;\;\;\;\;\;\;\;\;\;\; & \mbox{if $j_1 = j_2$ and $k_1 \neq k_2$} \\
							c^{m_1-1} & \mbox{if $j_1 \neq j_2$ and $k_1 = k_2$} \\
							c^n & \mbox{if $j_1 \neq j_2$ and $k_1 \neq k_2$} \\
							0 & \mbox{if $j_1 = j_2$ and $k_1 = k_2$}
						\end{array} \right\} m_1 = m_2
					\end{array} \right. \label{lkk} \end{align}\begin{align}
\lonek{1}{n}{j_1}{m_2}{n}{j_2}{k_2} &\equiv \left\{ \begin{array}{ll}
							a^{n-m_2+1} \;\;\;\;\;\;\;\;\;\;\;& \mbox{if $j_1 = j_2$} \\
							c^n & \mbox{if $j_1 \neq j_2$ and $j_1 \neq k_2$} \\
							c^{m_2-1} & \mbox{if $j_1 \neq j_2$ and $j_1 = k_2$} 
						\end{array} \right. \label{lonek} \\
\lkone{m_1}{n}{j_1}{k_1}{1}{n}{j_2} &\equiv \left\{ \begin{array}{ll}
							b \;\;\;\;\;\;\;\;\;\;\;\;\;\;\;\;\;\;\;\;\;\; & \mbox{if $j_1 = j_2$} \\
							c^n & \mbox{if $j_1 \neq j_2$ and $j_2 \neq k_1$} \\
							c^{m_1-1} & \mbox{if $j_1 \neq j_2$ and $j_2 = k_1$} 
						\end{array} \right. \label{lkone}\\
\loneone{1}{n}{j_1}{1}{n}{j_2} &\equiv \left\{ \begin{array}{ll}
							c^n \;\;\;\;\;\;\;\;\;\;\;\;\;\;\;\;\;\;\;\; & \mbox{if $j_1 \neq j_2$} \\
							0 & \mbox{if $j_1 = j_2$.}
						\end{array} \right. \label{loneone}
\end{align}
Additionally, to penalize initial state uncertainty, for $j,k \in \mathcal{S}$, the following costs are adopted:
\begin{eqnarray} \label{isurCost}
I(j,k) &\equiv& \left\{ \begin{array}{ll}
				\tau \;\;\;\;\;\;\;\;\;\;\;\; & \mbox{ if $j \neq k$ } \\
				0 & \mbox{ if $j = k$. }
				\end{array} \right.			
\end{eqnarray}
In (\ref{lkk})-(\ref{loneone}), the constant $a>1$ represents the base of the exponentially increasing cost of delay when the change is in the correct direction from initial state, and $c>1$ represents the base of the exponential cost of incorrect detection. Parameter $b > 0$ is the fixed cost of either (i) early correct detection ({\em  false alarm}), or (ii)  incorrect final state. 
In (\ref{isurCost}), the parameter $\tau >0$ serves as the fixed cost of initial state uncertainty. 
%This also includes 
%the factor $D/M(n) = 1/((n-1)(D-1)+1)$ .
%This factor creates 
%factor $1/M(n) = 1/D((n-1)(D-1)+1)$ common to all risk terms.
%Since common factors have no effect on the determination of the minimum risk hypothesis  in (\ref{minrisk}), 
%%the common factor of $1/M(n)$ 
%it is omitted without loss of generality.
%in the calculation of all risks for computational efficiency.
 
Applying  costs (\ref{lonek}) and (\ref{loneone}) to (\ref{R1n}), Bayes' rule (\ref{bayes}), (\ref{initposterior0}), as well as equally likely priors, the risk associated with choosing hypothesis $\hypone{1}{n}{j}$ at time $n$ is given by
\begin{align} 
 \riskone{1}{n}{j} 
& = \;\;  \frac{\pi}{P(X_{1,n} = x_{1,n})} \left[  \sum_{ \{ s  \in \mathcal{S}_{-j} \}} \sum_{i=2}^n a^{n-i+1}   P(X_{1,n} = x_{1,n} | \hypk{i}{n}{j}{s}) \right. \nonumber \\
& \hspace{20pt} + \sum_{ \{ r \in \mathcal{S}_{-j} \} } \Bigg{(}   c^n P( X_{1,n} = x_{1,n} | \hypone{1}{n}{r})   +  \sum_{ \{ s  \in \mathcal{S} | s_{-j-r}  \} } \sum_{i=2}^{n} c^n P(X_{1,n} = x_{1,n} | \hypk{i}{n}{r}{s}) \nonumber \\
& \hspace{75pt}+ \left.  \sum_{i=2}^n c^{i-1} P(X_{1,n} = x_{1,n} | \hypk{i}{n}{r}{j})    \Bigg{)} \right].  \label{recursive-R1n}
\end{align}
Using the indepedance of $X_{1,n+1}$,
it can be shown that  $\riskone{1}{n+1}{j}$ can be recursively updated from $\riskone{1}{n}{j}$ \cite{Falt17}. 
% by using recursions provided in Appendix~B. 
%and can be updated from $\riskone{1}{n}{j}$ using recursions provided in Appendix~A.
Similarly, the minimum risks over $\riskk{l}{n}{j}{k} , 1 < l \leq n, (j,k) \in \setsquare$ can also be tracked recursively. 
%The method is to recursively track the minimum of the risk associated with the current minimum risk change time and that associated with change time $n$. 
Let $m$ represent the minimum risk change time at time $n$ and let $m'$ be the minimum risk change time corresponding to time $n-1$. The minimum risk corresponding to change time $m$ is then updated at time $n$ according to
\begin{eqnarray}
\riskk{m}{n}{j}{k} = \min \{ \; \riskk{n}{n}{j}{k}, \; \riskk{{m'}}{n}{j}{k} \; \}. \label{minRnnRkpn}
\end{eqnarray}
In  (\ref{minRnnRkpn}), the Bayes risk $\riskk{m}{n}{j}{k}$ can be expressed using using costs (\ref{lkk}), (\ref{lkone}), and (\ref{isurCost}), 
%$\riskk{n}{n}{j}{k}$ can be expressed as
%\begin{align} 
%& \riskk{n}{n}{j}{k} \nonumber \\
%  = \;\; & \frac{\pi}{P(X_{1,n} = x_{1,n})} \Bigg{[}  \sum_{i=2}^{n-1} \Big{(} a^{n-i+1}P(X_{1,n} = x_{1,n} | \hypk{i}{n}{j}{k}) \Big{)} + b P(X_{1,n} = x_{1,n}| \hypone{1}{n}{j}) \nonumber \\
%& + \sum_{ \{ s \in \mathcal{S}_{-j-k} \} } \Bigg{(} \sum_{i=2}^{n} \Big{(} b P(X_{1,n} = x_{1,n} | \hypk{i}{n}{j}{s}) \Big{)}   \Bigg{)}   \nonumber \\
%& + \sum_{ \{ r \in \mathcal{S}_{-j-k}  \} } \Bigg{(} c^{n} P(X_{1,n} = x_{1,n} | \hypone{1}{n}{r}) \Bigg{)} + c^{n-1} P(X_{1,n} = x_{1,n} | \hypone{1}{n}{k}) \nonumber \\
%& + \sum_{ \{ r \in \mathcal{S} _{-j-k} \}  } \Bigg{(}  \sum_{i=2}^{n} \Big{(} c^{n-1} P(X_{1,n} = x_{1,n} | \hypk{i}{n}{r}{k}) \Big{)} + \sum_{i=2}^{n} \Big{(} c^{i} P(X_{1,n} = x_{1,n} | \hypk{i}{n}{r}{j}) \Big{)} \Bigg{)} \nonumber \\
%& + \sum_{i=2}^{n} \Big{(} c^i P(X_{1,n} = x_{1,n} | \hypk{i}{n}{k}{j})\Big{)} + \sum_{ \{ s \in \mathcal{S}_{-j-k}  \} } \Bigg{(} \sum_{i=2}^{n} \Big{(} c^n P(X_{1,n} = x_{1,n} | \hypk{i}{n}{k}{s}) \Big{)} \Bigg{)}   \nonumber \\
%& +  \sum_{ \{ r \in \mathcal{S}_{-j-k}  \} } \Bigg{(} \sum_{ \{ s \in \mathcal{S}_{-j-k}  \} } \Bigg{(} \sum_{i=2}^{n} \Big{(} c^n P(X_{1,n} = x_{1,n} | \hypk{i}{n}{r}{s}) \Big{)}    \Bigg{)} \Bigg{)} \Bigg{]} \nonumber \\
%& + t \frac{\sum_{ \{ r \in \mathcal{S}_{-j} \} } \prod_{i=1}^{n-1} f_r(x_i)}{\sum_{ \{ r \in \mathcal{S} \} } \prod_{i=1}^{n-1} f_r(x_i)} \label{recursive-Rnn}
%\end{align}
\begin{align}
& \riskk{m}{n}{j}{k} %\nonumber \\
= \;\;  \frac{\pi}{P(X_{1,n} = x_{1,n})} \Bigg{[} \sum_{i=2}^{m-1} a^{m-i} P(X_{1,n} = x_{1,n} | \hypk{i}{n}{j}{k}) \nonumber \\
& + \sum_{i=m+1}^{n} b P(X_1^n = x_1^n | \hypk{i}{n}{j}{k})  + b P(X_{1,n} = x_{1,n} | \hypone{1}{n}{j}) \nonumber \\
&  + \sum_{ \{ s \in \mathcal{S}_{-j-k} \}  } \Bigg{(} \sum_{i=2}^{m}  b P(X_{1,n} = x_{1,n} | \hypk{i}{n}{j}{s})   + \sum_{i=m+1}^{n} b P(X_{1,n} = x_{1,n} | \hypk{i}{n}{j}{s}) \Bigg{)} \nonumber \\
& + \sum_{ \{ r \in \mathcal{S}_{-j-k}  \} }  c^{n} P(X_{1,n} = x_{1,n} | \hypone{1}{n}{r}) + c^{m-1} P(X_{1,n} = x_{1,n} | \hypone{1}{n}{k}) \nonumber \\
& + \sum_{ \{ r \in \mathcal{S} _{-j-k} \}  } \Bigg{(}  \sum_{i=2}^{m} c^{m-1} P(X_{1,n} = x_{1,n} | \hypk{i}{n}{r}{k})   + \sum_{i=m+1}^{n} c^{i-1} P(X_{1,n} = x_{1,n} | \hypk{i}{n}{r}{k}) \nonumber \\
& \hspace{10pt}+  \sum_{i=2}^{m} c^{n-m+1}c^{i-1} P(X_{1,n} = x_{1,n} | \hypk{i}{n}{r}{j}) + \sum_{i=m+1}^{n} c^n P(X_{1,n} = x_{1,n} | \hypk{i}{n}{r}{j}) \nonumber \\
& \hspace{10pt}+ \sum_{i=2}^{m} c^{n-m+1}c^{i-1} P(X_{1,n} = x_{1,n} | \hypk{i}{n}{k}{j}) + \sum_{i=m+1}^{n} c^{n-i+1}c^{m-1} P(X_{1,n} = x_{1,n} | \hypk{i}{n}{k}{j})  \Bigg{)} \nonumber \\
& + \sum_{ \{ s \in \mathcal{S}_{-j-k}  \} } \Bigg{(} \sum_{i=2}^{m} c^n P(X_{1,n} = x_{1,n} | \hypk{i}{n}{k}{s}) + \sum_{i=m+1}^{n} c^{n-i+1}c^{m-1} P(X_{1,n} = x_{1,n} | \hypk{i}{n}{k}{s})   \Bigg{)} \nonumber \\
& + \sum_{ \{ r \in \mathcal{S}_{-j-k}  \}  }  \sum_{ \{ s \in \mathcal{S}_{-j-k}  \} } \sum_{i=2}^{n} c^n P(X_{1,n} = x_{1,n} | \hypk{i}{n}{r}{s})  \Bigg{]} 
 + \tau \frac{\sum_{ \{ r \in \mathcal{S}_{-j} \}  } \prod_{i=1}^{m-1} f_r(x_i)}{\sum_{ \{ r \in \mathcal{S} \} } \prod_{i=1}^{m-1} f_r(x_i)} \label{recursive-Rkp}
\end{align}
and are updated recursively using the most recent observation $x_{n}$.  For $\riskk{n}{n}{j}{k}$ in (\ref{minRnnRkpn}), the sums from $m+1$ to $n$ in (\ref{recursive-Rkp}) are omitted.
The detailed recursions can be found in  \cite{Falt17}.  
In short,  the minimum-risk hypothesis among the $M(n)$ possible change scenarios can be tracked recursively at time $n$ with constant time complexity by tracking the minimum of $\riskone{1}{n}{j}, j \in {\cal S}$ and  $\riskk{m}{n}{j}{k}, (j,k) \in \setsquare$, a total of $D + 2 \frac{D!}{(D-2)!}$ of the $D + (n-1)\frac{D!}{(D-2)!}$ risks which comprise the state of the detector. While recursive tracking of Bayes risk is possible for more general cost structures,  (\ref{lkk})-(\ref{isurCost}) is shown next to possess desirable properties. 

\section{Cost Parameter Bounds and Tradeoffs} \label{sec-design}

\subsection{Expected Risks}

Recall that cost increases exponentially with the number of observations from the true change time $m$, where the base of the exponents are $a$ and $c$. We establish bounds on parameter values  to ensure convergence of risks corresponding to hypotheses with the correct initial state and divergence of risks corresponding to hypotheses with incorrect initial states. From the risk expressions (\ref{Rkn}),  (\ref{R1n}),  (\ref{recursive-R1n}), (\ref{recursive-Rkp}) and Bayes rule, it can be noted that all risk terms excluding the initial state uncertainty risk terms have the common factor $1/P(X_{1,n} = x_{1,n})$. To establish the parameter bounds, we proceed by finding the expectation of all risks with the $1/P(X_{1,n} = x_{1,n})$ factor removed and  identify conditions for their convergence or divergence. After that, we show that these conditions extend to expectations of the risk itself.

Using conditional independence, taking expectations of the likelihoods (\ref{posterior_one}) and (\ref{posterior_k}), for $1 < m_1  \leq n$, $1 < m_2 \leq n $, and $(j,k), \; (r,s) \in \setsquare$ we obtain  
\begin{align}
& \mathbb{E}_{\hypk{m_2}{n}{j}{k}} [ P( X_{1,n} | \hypk{m_1}{n}{r}{s})  ] \nonumber \\
= & \left\{ \begin{array}{ll} 
\prod_{i = 1}^{m_1-1} \mathbb{E}_j [ f_r(X_i) ] \prod_{i=m_1}^{m_2-1} \mathbb{E}_j [ f_s(X_i) ] \prod_{i=m_2}^{n} \mathbb{E}_k [ f_s(X_i) ]   \quad & \mbox{for $m_1<m_2$,} \\
\prod_{i=1}^{m_2-1} \mathbb{E}_j [ f_r(X_i) ] \prod_{i=m_2}^{m_1-1} \mathbb{E}_k [ f_r(X_i) ] \prod_{i=m_1}^{n} \mathbb{E}_k [ f_s(X_i) ]   \quad & \mbox{for $m_1>m_2$, and} \\
\prod_{i=1}^{m_1-1} \mathbb{E}_j [ f_r(X_i) ]  \prod_{i=m_1}^{n} \mathbb{E}_k [ f_s(X_i) ]   \quad & \mbox{for $m_1=m_2$,} 
\end{array} \right.  \label{exppostkk}
\end{align}
where $\mathbb{E}_i[ \cdot ], \; i \in \mathcal{S}$ denotes expectation with respect to the distribution $f_i$, and $\mathbb{E}_{H1}[\cdot]$ denotes expectation with respect to the hypothesis $H1$. 
Similarly, we get
\begin{align}
\mathbb{E}_{\hypone{1}{n}{j}} [ P( X_{1,n} | \hypk{m_1}{n}{r}{s}) ] & = \prod_{i=1}^{m_1-1} \mathbb{E}_j [f_r(X_i)] \prod_{i=m_1}^{n} \mathbb{E}_j [ f_s(X_i)], \label{exppostk1}\\
\mathbb{E}_{\hypk{m_2}{n}{j}{k}} [ P( X_{1,n} | \hypone{1}{n}{r}) ] & = \prod_{i=1}^{m_2-1} \mathbb{E}_j [f_r(X_i)] \prod_{i=m_2}^{n} \mathbb{E}_k [ f_k(X_i)], \mbox{ and} \label{exppost1k} \\
\mathbb{E}_{\hypone{1}{n}{j}} [ P( X_{1,n} | \hypone{1}{n}{r})  ] & = \prod_{i=1}^n \mathbb{E}_j [f_r(X_i)]. \label{exppost11}
\end{align}
We rewrite (\ref{exppostkk})-(\ref{exppost11}) as
\begin{align}
& \mathbb{E}_{\hypk{m_2}{n}{j}{k}} [P(  X_{1,n} | \hypk{m_1}{n}{r}{s}) ] \nonumber \\
= & W_{m_2,n}(j,k) \left\{ \begin{array}{ll}
\dval{j}{r}{j}^{m_1-1} \dval{j}{s}{j}^{m_2-m_1} \dval{k}{s}{k}^{n-m_2+1} \quad & \mbox{for $m_1<m_2$} \\
\dval{j}{r}{j}^{m_2-1} \dval{k}{r}{k}^{m_1-m_2} \dval{k}{s}{k}^{n-m_1+1}  & \mbox{for $m_1>m_2$} \\
\dval{j}{r}{j}^{m_1-1} \dval{k}{s}{k}^{n-m_1+1}  & \mbox{for $m_1=m_2,$}
\end{array} \right. \label{exphypkk}
\end{align}
\begin{align}
\mathbb{E}_{ \hypone{1}{n}{j}} [P( X_{1,n}| \hypk{m_1}{n}{r}{s} )  ] & = W_{1,n}(j) \dval{j}{r}{j}^{m_1-1} \dval{j}{s}{j}^{n-m_1+1}, \\
\mathbb{E}_{\hypk{m_2}{n}{j}{k}} [P( X_{1,n}|\hypone{1}{n}{r})  ] & = W_{m_2,n}(j,k) \dval{j}{r}{j}^{m_2-1} \dval{k}{r}{k}^{n-m_2+1}, \\
\mathbb{E}_{\hypone{1}{n}{j}} [P( X_{1,n}|\hypone{1}{n}{r} ) ] & = W_{1,n}(j) \dval{j}{r}{j}^n \label{exphyponeone}
\end{align}
where constants 
\begin{align*}
%W_{1,n}(j) & \stackrel{\Delta}{=} \mathbb{E}_{\hypone{1}{n}{j}} [ P(X_{1,n} | \hypone{1}{n}{j}  ) ], \\
W_{1,n}(j) & \stackrel{\Delta}{=}  \mathbb{E}_j [f_j (X)]^n,  \\
%W_{m_2,n}(j,k) & \stackrel{\Delta}{=} \mathbb{E}_{\hypk{m_2}{n}{j}{k}} [  P( X_{1,n}| \hypk{m_2}{n}{j}{k}  ) ], \; \mbox{and} \\
W_{m_2,n}(j,k) & \stackrel{\Delta}{=} \mathbb{E}_j [f_j (X)]^{m_2-1} \mathbb{E}_k [f_k(X)]^{n - m_2+1}, \; \mbox{and}
\end{align*} 
\begin{align}
\dval{i}{j}{k} &  \stackrel{\Delta}{=}  \frac{\mathbb{E}_i \left[ f_j(X) \right]}{\mathbb{E}_i \left[ f_k(X) \right]} \quad \quad \quad 	 i,j,k \in \mathcal{S}. \label{d}
\end{align}
We note in Eq.~(\ref{d}), 
\begin{lemma} \label{lemma1}
For any set of distinct pdfs $\{ f_j | j \in \mathcal{S} \}$ such that $\langle f_j,f_j \rangle = \langle f_k,f_k \rangle \; \forall j,k \in \mathcal{S}$, $\dval{j}{j}{k} > 1$ and  $\dval{j}{k}{j} < 1 $ for all $(j,k) \in \setsquare $.
\end{lemma}

\noindent
Proof: See Appendix~A.

%Though functions of the length of the test, $n$, $W_{1,n}(j)$ and $W_{m,n}(j,k)$, only depend on the known PDFs $f_j$ and $f_k$ and assumed change time $m$ and serve as common factors to the different actual outcomes, $X_{1,n} = x_{1,n}$, unlike $\dval{j}{r}{j}$, whose numerator depends on the $r$th observed sample's PDF which may have different change time, initial PDF and final PDF.

The $W_{1,n}(j,k)$ and $W_{m,n}(j,k)$ are functions of the length of the test, $n$, and only depend on the known PDFs $f_j $ and $f_k$ and the assumed change time $m$. In the following, these functions serve as common factors when calculating the expected value of the conditional risks. In contrast, the values $d_i(j,k)$ serve to represent the factor by which the expected value of different likelihoods (\ref{exppostkk})-(\ref{exppost11}) differ from the functions $W_{1,n}(j,k)$ and $W_{m,n}(j,k)$.

Multiplying (\ref{recursive-Rkp}) by $P(X_{1,n})$, averaging over $X_{1,n} \sim \hypk{m}{n}{j}{k}$ for $1 < m \leq n$, 
and applying (\ref{exphypkk})-(\ref{exphyponeone}) 
yields the likelihood-weighted expected risk for the {\em correct detection} hypothesis
\begin{align}
& \mathbb{E}_{\hypk{m}{n}{j}{k}}[ \riskk{m}{n}{j}{k} P(X_{1,n} ) ] \nonumber \\
 = & \; \pi W_{m,n}(j,k) \Bigg{[} \sum_{i=2}^{m-1} (a \dval{j}{k}{j})^{m-i} + \sum_{i=m+1}^{n}  b (\dval{k}{j}{k})^{i-m} + b (\dval{k}{j}{k})^{n-m+1} \nonumber \\
& + \sum_{\{ s \in \mathcal{S}_{-j-k} \}} \Big{(} \sum_{i=2}^{m} b \dval{j}{s}{j}^{m-i} \dval{k}{s}{k}^{n-m+1} + \sum_{i=m+1}^{n} b \dval{k}{j}{k}^{i-m} \dval{k}{s}{k}^{n-i+1} \Big{)} \nonumber 
\end{align}
\begin{align}
& + \sum_{ \{ r \in \mathcal{S} _{-j-k} \} } c^n \dval{j}{r}{j}^{m-1} \dval{k}{s}{k}^{n-m+1} + c^{m-1} \dval{j}{k}{j}^{m-1} \nonumber \\
& + \sum_{ \{ r \in \mathcal{S} _{-j-k} \}  } \Bigg{(} \sum_{i=2}^{m} c^{m-1} \dval{j}{r}{j}^{i-1} \dval{j}{k}{j}^{m-i}  + \sum_{i=m+1}^{n} c^{i-1} \dval{j}{r}{j}^{m-1} \dval{k}{r}{k}^{i-m}  \nonumber \\
& + \sum_{i=2}^{m} (c \dval{k}{j}{k})^{n-m+1} (c \dval{j}{k}{j})^{i-1} + \sum_{i=m+1}^{n} c^n \dval{r}{j}{r}^{m-1} \dval{r}{k}{r}^{i-m} \dval{k}{j}{k}^{n-i+1}  \nonumber \\
& + \sum_{i=2}^{m} (c \dval{k}{j}{k})^{n-m+1} (c \dval{j}{r}{j})^{i-1}  + \sum_{i=m+1}^{n} (c \dval{k}{j}{k})^{n-i+1} (c \dval{j}{k}{j})^{m-1} \Bigg{)} \nonumber \\
& + \sum_{\{ s \in \mathcal{S}_{-j-k} \}} \Bigg{(} \sum_{i=2}^{m} c^n \dval{j}{k}{j}^{i-1} \dval{j}{s}{j}^{m-i} \dval{k}{s}{k}^{n-m+1} \nonumber \\
& \hspace{100pt} + \sum_{i=m+1}^{n} (c \dval{j}{k}{j})^{m-1} (c \dval{k}{s}{k})^{n-i+1} \Bigg{)} \nonumber \\
& + \sum_{\{ r \in \mathcal{S} _{-j-k} \} } \sum_{\{ s \in \mathcal{S}_{-j-k} \}} \Bigg{(} \sum_{i=2}^m c^n \dval{j}{r}{j}^{i-1} \dval{j}{s}{j}^{m-i} \dval{k}{s}{k}^{n-m+1} \nonumber \\
& \hspace{100pt} + \sum_{i=m+1}^{n} c^n \dval{j}{r}{j}^{m-1} \dval{k}{r}{k}^{i-m} \dval{k}{s}{k}^{n-i+1}  \Bigg{)}  \Bigg{]} \nonumber \\
& + \mathbb{E} \left[ P(X_{1,n} ) \tau \frac{\sum_{ \{ r \in \mathcal{S}_{-j} \}} \prod_{i=1}^{m-1} f_r(X_i)}{\sum_{ \{ r \in \mathcal{S} \} } \prod_{i=1}^{m-1} f_r(X_i)}  \right]  \label{ERkn}
\end{align}
which is comprised of geometric series terms, each a function of cost parameters and PDFs. The initial state uncertainty risk in the final term of (\ref{ERkn}) tends to zero as $m$ gets large, and thus has vanishing influence on test behaviour for large change times.

%For the expected risk of a hypothesis under the conditions of incorrect detection, we proceed by taking the expected value of $\riskk{m}{n}{p}{q} P(X_{1,n} = x_{1,n})$ conditioned on $\hypk{m}{n}{j}{k}$ being true, for $ (p,q), \; (j,k) \in \setsquare$ and $p \neq j$. This value can be determined using (\ref{Rkn}), the cost structure (\ref{lkk})-(\ref{loneone}) and (\ref{isurCost}), and (\ref{exphypkk})-(\ref{exphyponeone}). 

\subsection{Parameter Choices for Large Change Times}

To provide insight into the time sequential notions of {\em correct} and {\em incorrect detection}, we first focus on the {\em large change-time regime} where after observing for a long time, the change has yet to occur. That is, in this subsection, we assume that $n \rightarrow \infty$ while $n-m$ remains finite, %we ignore cases where a change had already occurred an asymptotically long time ago
and thereby avoid situations where an incorrect detection corresponds to a vanishingly small transient initial state, an assumption consistent  with the chosen cost structure. 
Observing (\ref{ERkn}), every term involving an exponential cost is a product of terms of the form $a \: \dval{r}{s}{r}$ or $c \: \dval{r}{s}{r}$ for $(r,s) \in \setsquare$.  If we choose parameters $a$ and $c$ to be such that each of these individual terms is less than one, then each of the geometric series terms with exponential costs will asymptotically converge, and finite shifts in change times in (\ref{ERkn}) do not affect convergence. The terms which include the fixed cost, $b$, will converge to a finite value for any $0 < b < \infty$.  % with respect to $\pi W_{m,n}(j,k).$
Thus, under the large change-time regime, the expected risk  Eq.~(\ref{ERkn}) can be shown to converge by choosing $b <  \infty$, $\tau < \infty$, $1 < a < d_{min}$, and $1< c < d_{min}$, where 
\begin{align}
d_{min} = \min_{\{ (r,s) \in \setsquare \}} \dval{r}{r}{s}. \label{dmin}
\end{align}
We note that $d_{min} > 1$ by Lemma \ref{lemma1} and we recall that $a>1$ and $c>1$ are assumed in the cost function definition. 
It can similarly be shown that choosing $1 < c < d_{min}$ will result in the expected risks for hypotheses with incorrect initial state diverging asymptotically by considering $\mathbb{E}_{\hypk{m}{n}{j}{k}} \left[ \riskk{m}{n}{p}{q} P(X_{1,n})\right]$ for $p \neq j$. The process of calculating $\mathbb{E}_{\hypk{m}{n}{j}{k}} \left[ \riskk{m}{n}{p}{q} P(X_{1,n})\right]$ for $p \neq j$ is the same as was used to find Eq. (\ref{ERkn}), but in this case several geometric series terms will diverge exponentially with base $c$ as $n$ increases. % Applying Eqs. (\ref{cov1})-(\ref{cov2}) to $\mathbb{E}_{\hypk{m}{n}{j}{k}} \left[ \riskk{m}{n}{p}{q} P(X_{1,n})\right]$, we conclude that under the large change-time regime $\mathbb{E}_{\hypk{m}{n}{j}{k}} \left[ \riskk{m}{n}{p}{q} \right]$ diverges asymptotically for $p \neq j$.

%In $ \riskk{m}{n}{p}{q} P(X_{1,n} )$, taking the expectation of the term $ \pi \lkk{m}{n}{p}{q}{m}{n}{j}{k} P(  X_{1,n} | \hypk{m}{n}{j}{k} )$ over $ \hypk{m}{n}{j}{k}$ results in at least $ \pi W_{m,n}(j,k) c^{m-1}$, where the cost $c^{m-1}$ is the smallest cost which $\lkk{m}{n}{p}{q}{m}{n}{j}{k}$ may assume for $p \neq j$. 
%As such, when calculating $\mathbb{E}_{\hypk{m}{n}{j}{k}} \left[ \riskk{m}{n}{p}{q} P(X_{1,n}) \right]$ similarly to (\ref{ERkn}), the geometric series terms will diverge exponentially at a rate which is at least exponential with base $c$ as $n$ grows. Substituting $\riskk{m}{n}{p}{q}$ for $\riskk{m}{n}{j}{k}$ in (\ref{cov2}) and taking expectation over $\hypk{m}{n}{j}{k}$ being true, we can conclude that $\mathbb{E}_{\hypk{m}{n}{j}{k}} \left[ \riskk{m}{n}{p}{q} \right]$ will diverge as $n$ increases under the large change-time regime.
%It can similarly be shown that, under the large change-time regime, all expected risks for hypotheses which correspond to the wrong initial state will diverge asymptotically with respect to $\pi W_{m,n}(j,k)$.

Next, we need to establish
\begin{lemma} \label{lemma2}
Assume that each $f_j, j \in \cal{S}$ has finite variance. For $(j,k) \in \setsquare$,  if the expected value of the likelihood-weighted risk $\mathbb{E} \big{[} \riskk{m}{n}{j}{k} P(X_{1,n}) \big{]}$  converges or diverges as $n$ gets large, $\mathbb{E}[ \riskk{m}{n}{j}{k}]$ will also converge or diverge, respectively.
\end{lemma}
\noindent
Proof: See Appendix~B.

We then  have   
\begin{theorem} \label{theorem1}
Under the conditions 
\begin{equation} \label{cond1}
1 < a <  d_{min} \;\mbox{and} \; 1 < c <  d_{min},
\end{equation}
 $0 < b < \infty$, and $0 < \tau < \infty$, the probability of incorrect detection of the proposed procedure (\ref{minrisk}) converges to zero asymptotically for the large change-time regime as the number of observations gets large.
\end{theorem}

\begin{proof}
According to Lemma~\ref{lemma2}, for establishing convergence or divergence, we may consider the average risks of the form (\ref{ERkn}). Applying the conditions (\ref{cond1}) to (\ref{ERkn}), at least one of the $(n-1)$ expected risks corresponding to correct detection converges to a steady state. This means that the minimum Bayes risk of correct detection in (\ref{minrisk})  must be bounded, and denote its value by $r_{\min} < \infty$.  Consider the random variable $(r_{\min}  - \riskk{m}{n}{p}{q})_+$, where $\riskk{m}{n}{p}{q}$ is the risk associated with choosing $\hypk{m}{n}{p}{q}$, a hypothesis with incorrect initial state, and $(\cdot)_+$ denotes $\max(\cdot, 0)$. It follows from (\ref{cond1}) that $\riskk{m}{n}{p}{q}$ is a diverging incorrect detection hypothesis. Applying Markov's inequality,
\begin{equation}
P\left(  (r_{\min}  - \riskk{m}{n}{p}{q})_+ > \epsilon \right) < \frac{\mathbb{E} \left[ (r_{\min} - \riskk{m}{n}{p}{q})_+  \right]}{\epsilon}
\end{equation}
for any $\epsilon >0$.
Since $\mathbb{E} [\riskk{m}{n}{p}{q}] \rightarrow \infty$ under the large change-time regime, the expectation on the right hand side vanishes as $n\rightarrow \infty$ for any finite $r_{\min}$. As any positive probability on the left-hand side represents the probability that $\riskk{m}{n}{p}{q}$ is the minimum risk hypothesis, we see that the incorrect detection probability vanishes to zero.
\end{proof}

\begin{remark}
Risks for hypotheses with incorrect initial state will stay large under the large change-time regime after a change occurs. In fact, a larger value of $c$ causes incorrect-detection risks to diverge more quickly, reducing the probability of incorrect detection. In addition, the time spent in an initial state where  incorrect-detection risks are small enough to incur an incorrect detection decrease as $c$ increases. 
\end{remark}

%%%%%%%%%%%%%%%%%%%%%%%%%% THEOREM 2 

\subsection{Parameter Choices for Initial State Uncertainty due to Small Change Times}

Theorem \ref{theorem1} is conditioned on the large change-time regime, i.e., each of the $\frac{D!}{(D-2)!}$ recursively tracked minimum-risk change times growing with $n$ before a change occurs, and is presented to justify the exponential cost structure. It is worth noting, however, that the assumption of a large change-time regime is not valid in all situations. If the parameters $a$ and $c$ are selected according to (\ref{cond1}), then the recursively tracked change times for risks corresponding to hypotheses with incorrect initial state should stay small since the algorithm identifies the minimum-risk hypothesis corresponding to change in each direction. Consequently, if the recursively tracked change time is small, then the asymptotic analysis leading to Theorem \ref{theorem1} does not apply. Additionally, for small change times the exponential cost structure does not always associate a large cost with incorrect detection, which is obviously undesirable. This motivates the use of the initial state uncertainty cost, $\tau$, which provides a mechanism to address the inherently large probability of incorrect detection for small change times.

For any particular risk corresponding to change, when the recursively tracked time is small, only a small number of samples can be used to determine the initial state of the sequence $X_{1,n}$. The cost parameter $\tau$ mitigates risk caused by early change-time hypotheses to prevent incorrect detections caused by this initial state uncertainty. As the number of initial state samples increases, the risk vanishes for correct-sided risks corresponding to change and remains large for incorrect-sided risks. Before we proceed, we define $f_{\mathcal{S}_{-j}}(X) = \frac{1}{D-1}\sum_{\{ r \in \mathcal{S}| r \neq j \}} f_r(X)$ 
%as the uniform linear combination of PDFs  $\left\{ r \in \mathcal{S} | r \neq j \right\}$, 
and $D_{KL}( f || g ) = \mathbb{E}_f \left[ \log \left( \frac{f(X)}{g(X)}\right) \right]$ as the Kullback-Leibler divergence between PDFs $f(X)$ and $g(X)$. The initial state uncertainty risk associated with the proposed procedure, (\ref{minrisk}), is shown to have the following property:

\begin{theorem} \label{theorem2}

Under the conditions (\ref{cond1}) and $0<b< \infty$, the following properties hold:

\noindent (i) the probability of incorrect detection of the proposed test can be made arbitrarily small by using a sufficiently large threshold $\tau$, 

\noindent (ii) for change times 
\begin{align}
m \geq \left \lceil{1 + \frac{ \log \left( \frac{\tau}{C_m} - 1 \right) + \log \left( D-1 \right)}{ D_{KL} \left( f_j || f_{\mathcal{S}_{-j}} \right) }} \right \rceil \label{logt} 
\end{align}
for some $j \in \mathcal{S}$, where $C_m$ is a finite constant, initial state uncertainty increases the expected correct detection delay as  $O({\log(\tau)})$,  and 
 
\noindent (iii) as
% each of the $D$ distributions become more similar, 
\begin{align}
\max_{\{ (j,k) \in \setsquare\}} D_{KL}(f_j||f_k) \rightarrow 0, \label{dist_similar}
\end{align}
initial state uncertainty increases expected delay as $O(\log(d_{min})^{-1})$,
where $d_{min}$ is given by (\ref{dmin}). 

\end{theorem}

\begin{proof}
See Appendix~B.
\end{proof}

\begin{remark}
The choice of parameter $\tau$ governs the trade off of incorrect detection probability against the ability to detect changes occurring prior to some minimum detectable change time.
\end{remark}

In Theorem \ref{theorem2}, a value of the minimum change time (\ref{logt}) is used to showcase performance trade-offs; however, the  value of this minimum change time is not useful 
for practical test design purposes since it is simply the average number of samples from the initial state required for the initial state uncertainty risk to vanish among other risk terms, which is expressed by Eq. (\ref{partonecondition}) in Appendix B. Rather, it is more useful to find the smallest change time $m$ satisfying
%\begin{align}
%P_{\hypone{1}{m-1}{j}} \left( \log \left( \frac{f_{j}^{m-1}(X_{1,m-1} )}{ f_{\mathcal{S}_{-j}}^{m-1} (X_{1,m-1}) } \right) > \log \left( \frac{\tau}{C_m}-1 \right) + \log \left( D-1 \right) \right) > 1 - \alpha, \label{prob_delay1}
%\end{align}
%where $0 \leq \alpha \leq 1$ is chosen. Equivalently,  
%$X_{1,m-1}$ is IID, (\ref{prob_delay1}) can equivalently be written as
\begin{align}
P_{\hypone{1}{m-1}{j}} \left( \sum_{i=1}^{m-1} \log \left( \frac{f_{j}(X_i)}{ f_{\mathcal{S}_{-j}} (X_i) } \right) < \log \left( \frac{\tau}{C_m}-1 \right) + \log \left( D-1 \right) \right) & < \alpha \label{prob_delay2}
\end{align}
for some chosen $0 \leq \alpha \leq 1$, which can be obtained from Eq. (\ref{partoneconditioneqv}) in Appendix~B.
This corresponds to the minimum change time for initial state uncertainty to affect average detection delay with probability $\alpha$. 
%As $\alpha \rightarrow 0$, the increase in delay 
%due to initial state uncertainty 
%vanishes. 
Unfortunately, computation of  the minimum $m$ satisfying (\ref{prob_delay2}) in general is difficult as it requires the cumulative distribution function of $\log \left( \frac{f_{j}(X_i)}{ f_{\mathcal{S}_{-j}} (X_i) } \right) $ conditioned on $\hypone{1}{m-1}{j}$. To avoid this, in Section~\ref{subsec-analysis} the increase in average detection delay due to initial state uncertainty risk is quantified instead.

\subsection{Delay and False Alarm} 
\label{sec_design_delay_fa}

It is worth further investigating the correct-sided risks for finite $n$ under the conditions (\ref{cond1}).  
Interpreting Eq.~(\ref{ERkn}) and Eq.~(\ref{ER1n}), Appendix~B, it can be noted that increasing parameter $b$ for fixed $a$ and $c$ serves as a mechanism to lower both probability of false alarm and probability of incorrect detection.  Under correct detection hypotheses,  in view of the terms including cost parameter $b$ in Eq.~(\ref{ERkn}) and  $\dval{k}{j}{k} < 1$, false alarms are controlled by the added risk (which incurs added delay).  Therefore,  increasing  $b$  while using values of $a$ and $c$ that satisfy Theorem~\ref{theorem1} can be used to achieve additional performance tradeoffs of decreased detection error probabilities at the expense of delay. 

Nevertheless, from Theorem \ref{theorem2},  there exists an irreducible minimum change time required for correct-sided initial state uncertainty risk to vanish, and it can be expressed directly in terms of the cost parameters and distributions. Since $C_m$ in (\ref{logt}) is not a function of cost parameter $\tau$,  
%(\ref{logt}) shows that 
this minimum change time only increases logarithmically in $\tau$, which indicates insensitivity to delay penalty. Taking (\ref{upperboundmaxmin}), Appendix~B, into account, we conclude that increasing $\tau$ to reduce the probability of incorrect detection would incur only a modest effect on delay in correct detection.

\section{Simulations and Numerical Results} \label{sec-results}

%% HERE WE NEED TO SUMMARIZE PREVIOUS RESULTS BEFORE MOVING ON.

% We have shown that:

%% test behaves as designed (sim which tracks risks)

%% compares well with CUSUM when initial state is strongly established

%% despite having 4 parameters, test design for good performance is easily achieved by:

%%% choose a,c arbitrarily (somewhat)
%%% choose b to trade off delay/fa
%%% choose t sufficiently large to prevent incorrect detections due to initial state uncertainty. 

\subsection{Simulation Description and Parameter Selection} \label{subsec-simdescription}

To illustrate the performance of the proposed test and its finite sample properties, i.e.,  Theorem \ref{theorem2}, Monte Carlo simulations were performed when initial state uncertainty cost $\tau$ is changed. We use $D=2$, which is sufficiently general for this purpose. An example featuring $D>2$ can be found  in \cite{Falt17}. Here, consider detecting a change in the mean of Gaussian PDF, $f_0$ to  $f_1$, i.e.,
\begin{align} \begin{array}{ll}
f_0 &= \mathcal{N}(0, \sigma^2) \\
f_1 &= \mathcal{N}(\mu, \sigma^2)
\end{array} \label{PDF_test1}
\end{align}
where 
$\mu$ is the mean shift after the change in distribution and 
$\sigma^2$ is the variance. 
%The use of $D=2$ is sufficient for the purposes of illustrating the finite sample behaviour addressed by Theorem 2.
%, as Theorem 2 presents performance bounds concerned with the initial state uncertainty risk, and alternative initial states occur for any $D \geq 2$. 
Define signal-to-noise ratio as $\mbox{SNR} \equiv \mu^2 / \sigma^2$. From (\ref{PDF_test1}), (\ref{d}) and (\ref{dmin}) we obtain $d_{min} = \dval{0}{0}{1} = \dval{1}{1}{0} = e^{\mbox{SNR}/4}$. 
Using an SNR of $0\mbox{dB}$  yields $d_{min} = 1.2840$. To satisfy conditions (\ref{cond1}) for Theorem \ref{theorem1}, we choose parameter values $a = 1.05$ and $c = 1.25$ to associate a larger cost with incorrect detection than detection delay. A cost of false alarm $b = 10^{1.85}$ is chosen, using the trade-off between average detection delay and false alarm rate discussed in Section \ref{sec_design_delay_fa}, to achieve a false alarm rate of approximately 0.05 when the change time is $m=50$ and the initial state uncertainty cost is $\tau=0$. To illustrate how initial state uncertainty cost influences probability of incorrect detection of initial state and average detection delay, initial state uncertainty cost $\tau$ is varied from $10^1$ to $10^7$, and for each value of $\tau$, the test performance is simulated using $10^6$ Monte Carlo trials for changes times varying from 5 to 100 in increments of 5.  For each trial, the initial state of the sequence is chosen randomly, and assumes either $f_0$ or $f_1$ with equal probability. The results are presented in Section \ref{subsec-MCresults}. In Section \ref{subsec-analysis}, the performance bound on the probability of incorrect detection and the average detection delay presented in Theorem \ref{theorem2} is compared with observed results and trends. Finally, Section \ref{subsec-summaryresults} provides summary discussion.
%Figures \ref{plotchangetp1}, \ref{plotchangetp2}, and \ref{plotchangetp3} show the average detection delays, incorrect detection rates, and false alarm rates, respectively, obtained from the simulation. 

An obvious comparison point is to use a fixed sample size (FSS) hypothesis test (HT) to detect the initial state of the observed sequence. Once the initial state is identified, a traditional change detection scheme, i.e., one that assumes knowledge of the initial state, such as Page's CUSUM, could be used to assess detection delay. Since choosing the length of the FSS HT requires knowledge of the change time, which is unrealizable, this is only used as a benchmark; in Section \ref{subsec-MCresults}, the incorrect initial state detection probability of the FSS HT and the average detection delay and false alarm rate of CUSUM are presented alongside that of the proposed change detection scheme. For the FSS HT, the incorrect detection rate shown assumes a known change time and chosen to be the change time minus one.
%For reference, Figure \ref{plotchangetp1} also shows the average detection delay of CUSUM when the initial and final states are assumed to be known. 
For a meaningful  comparison, a  CUSUM threshold of  $10^{2.157}$ was chosen achieve approximately the same ARL to false alarm as the proposed scheme with the  above parameter values. %Figure \ref{plotchangetp2} also shows the incorrect detection rate of the FSS HT for determining the initial state assuming the change time is known, i.e. the length of the FSS HT is the change time minus one. Figure \ref{plotchangetp3} additionally shows the false alarm rate of CUSUM for the case when the initial and final states of the sequence are known. 

\subsection{Monte Carlo Simulation Results} \label{subsec-MCresults}

\begin{figure}[!tpb]
\centering
\includegraphics[scale=0.7]{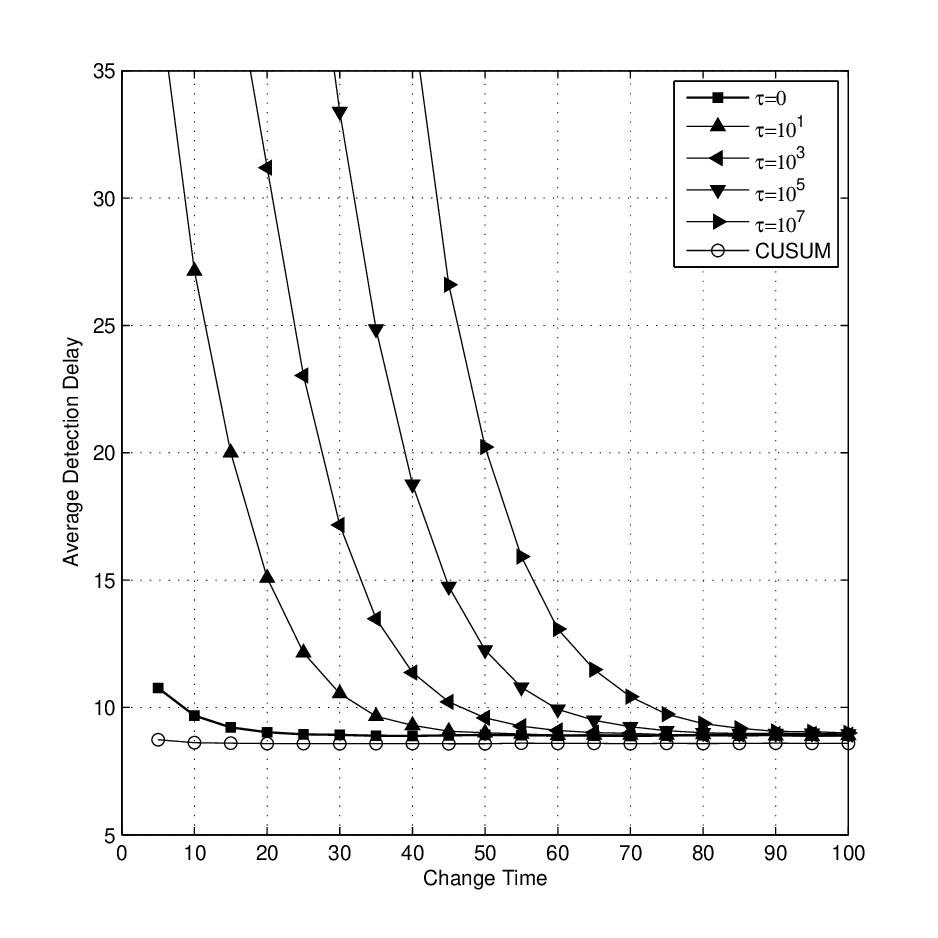}
\caption{Average detection delay versus change time of proposed test for different initial state uncertainty costs, $\tau$. The PDFs $f_0$ and $f_1$ in (\ref{PDF_test1}) indicate a change in the mean of a Gaussian distribution. Also shown is the average detection delay for CUSUM for the case where the initial and final states are assumed known. Parameter values are $a=1.05$, $c=1.25$, and $b=10^{1.85}$. CUSUM's threshold is $10^{2.157}$. }
\label{plotchangetp1}
\end{figure}

\begin{figure}[!tpb]
\centering
\includegraphics[scale=0.7]{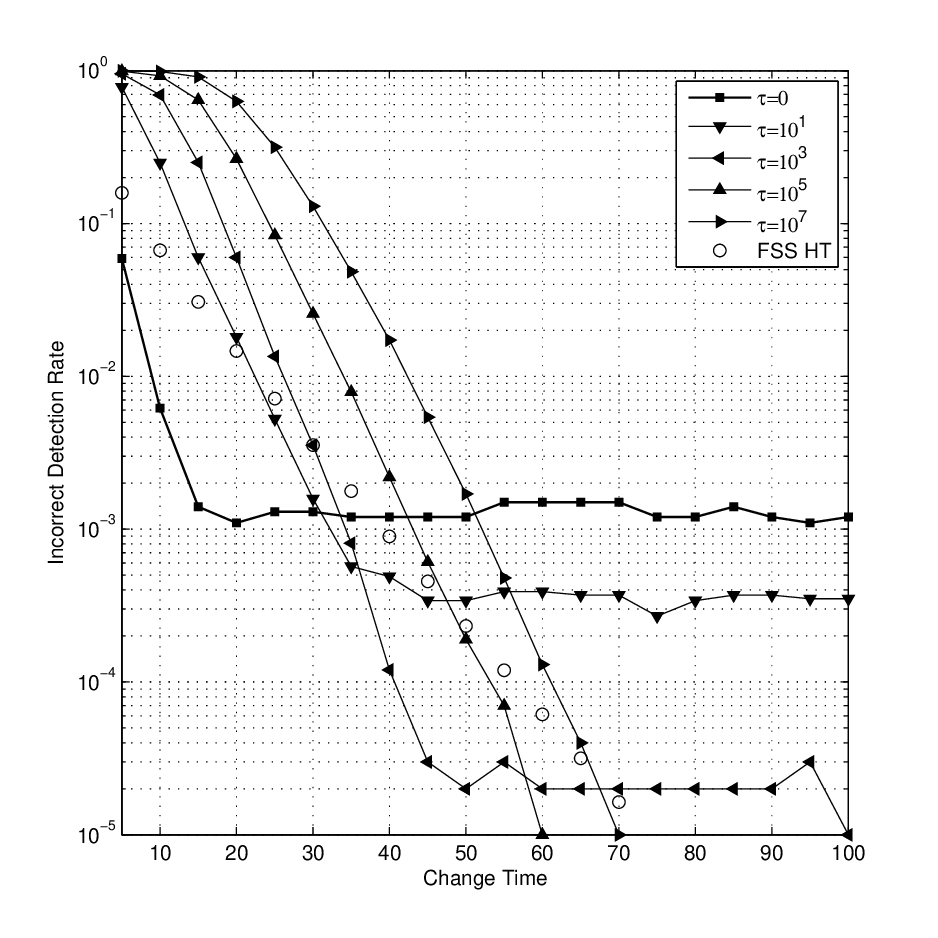}
\caption{Incorrect detection rate versus change time for different initial state uncertainty costs, $\tau$. The PDFs are those used in Fig.~\ref{plotchangetp1}.  Also shown is the incorrect detection probability of a fixed sample size hypothesis test of the initial distribution of a sequence assuming a known change time of $m$. Parameter values are $a=1.05$, $c=1.25$, and $b=10^{1.85}$. }
\label{plotchangetp2}
\end{figure}

\begin{figure}[!tpb]
\centering
\includegraphics[scale=0.7]{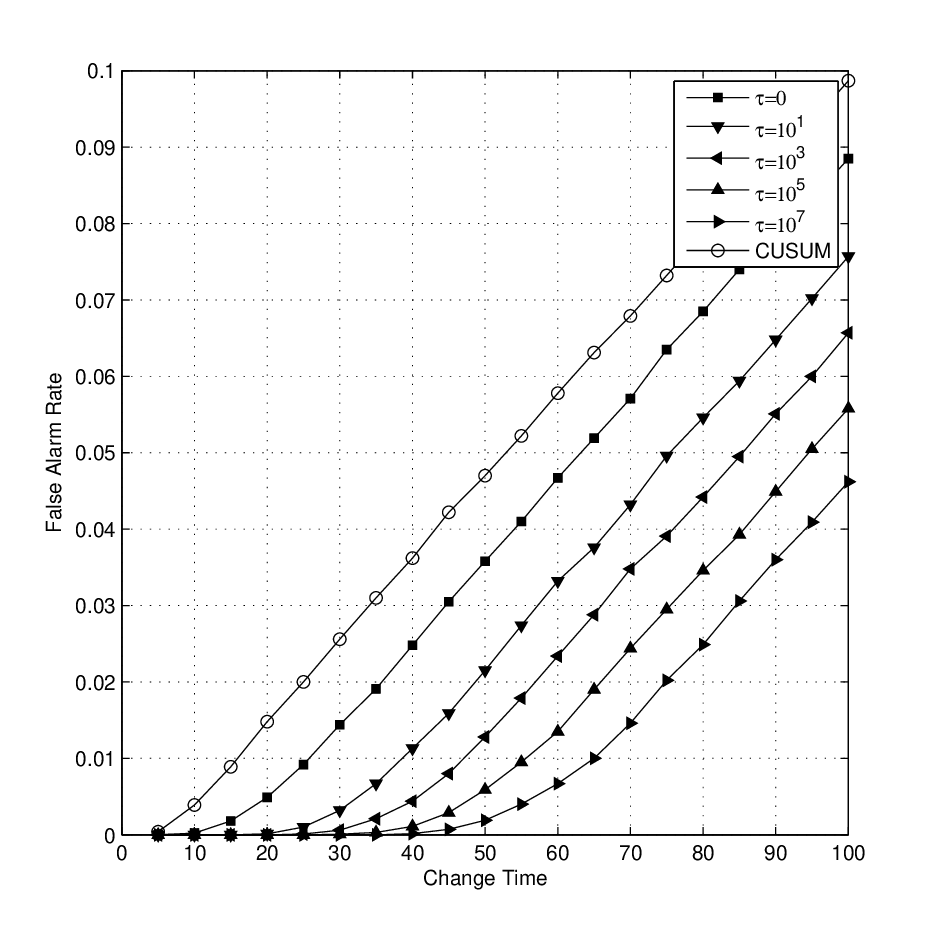}
\caption{False alarm rate versus change time for different initial state uncertainty costs, $\tau$. The PDFs are those used in Fig.~\ref{plotchangetp1}. Also shown in the false alarm rate of CUSUM, whose threshold is chosen such that CUSUM's ARL to false alarm is approximately that of the proposed change detection scheme. Parameter values are $a=1.05$, $c=1.25$, and $b=10^{1.85}$, and the threshold used for CUSUM is $10^{2.157}$.}
\label{plotchangetp3}
\end{figure}

For the PDFs (\ref{PDF_test1}) and parameter values selected in Section \ref{subsec-simdescription}, Figures \ref{plotchangetp1}, \ref{plotchangetp2}, and \ref{plotchangetp3} plot average detection delays, incorrect detection rates, and false alarm rates, respectively, achieved by the proposed change detection scheme in the Monte Carlo simulation.

From Figure 1 it is clear that average detection delay increases for any given change time when the initial state uncertainty cost, $\tau$,  is increased, which is consistent with Theorem \ref{theorem2}. When the change time is small, the average detection delay of the test for all $\tau>0$ is much larger than for $\tau=0$. However, as the change time increases, the average detection delay becomes insensitive to cost $\tau$ as expected. Additionally, it can be noted that the average detection delay of the proposed change detection scheme after an initial transient period is only slightly greater than that of CUSUM, which is known to be average-delay-optimal for known initial state and FAR constraint. The average detection delay of CUSUM is measured to be 8.586, while that of  proposed scheme once the initial state is established is 8.903, which is only 3.69 \% larger. 

Observing Figure 2, it is clear that each value of $\tau$ achieve a minimum probability of incorrect detection, which decreases as $\tau$ is increased. 
At the same time, increasing $\tau$ also increases the minimum change time for when this minimum incorrect detection rate occurs. These results are consistent with Theorem \ref{theorem2}. For example, for $\tau=10^3$, a minimum incorrect detection rate of approximately \num{2e-5} is achieved and it reaches this minimum at a change time of approximately sample 45. When $\tau=0$ is used, a much larger incorrect detection rate floor of approximately \num{1.3e-3} is achieved earlier, by sample 15. These results clearly illustrate the trade-off  between the probability of incorrect detection of initial state and the test's ability to detect early change times. Comparing the proposed change detection scheme's incorrect initial state detection rate to that of the FSS HT, it can be observed that at for each change time there is a  value of $\tau$ for the proposed test that achieves a lower incorrect detection rate than the FSS HT. Not surprisingly, no value of $\tau$ is universally better than the FSS HT which requires initial state knowledge. 
%For the purposes of test design, it would thus be beneficial to have some knowledge of the change time to select the parameter value of $\tau$, although when none is available, choosing $\tau=$0$ will ensure the best performance with regards to detecting small change times.%Additionally, it is worth noting that, for $\tau=$10^1$ and $\tau=$10^3$, the change time where the incorrect detection floor in Figure \ref{plotchangetp2} is reached coincides approximately with the change time at which the increase in average detection delay compared to the delay for $\tau=$0$ becomes small for in Figure \ref{plotchangetp1}. 

Figure 3 shows the false alarm rates achieved by the proposed change detection scheme for different values of $\tau$. It is clear that after an initial transient period during which the initial state is established, false alarm rate increases linearly with change time. This is consistent with Theorem 2, that specifies how the initial transient period increases with $\tau$. It can also be noted that, following the initial transient behaviour, the rate at which false alarm rate increases with change time is approximately invariant with $\tau$, which indicates that the parameter $\tau$ has a negligible effect on the ARL to false alarm. Recalling that an increase in average detection delay caused by initial state uncertainty risk vanishes following initial transient behaviour of the test, we conclude that the initial state uncertainty risk does not affect the trade-off between the probability of false alarm and average detection delay.
%if the change time is large enough for the test to have a high probability of correctly establishing the initial state, increasing the change time will increase the false alarm rate linearly. 
Finally, it can be noted that the rate at which the false alarm rates increase for both CUSUM and the proposed change detector are approximately the same. This was intended, as it was desired for both tests to have approximately the same ARL to false alarm for benchmarking purposes. It should be noted that Figure 3 presents the false alarm rate of CUSUM under the assumption that the initial and final states are known. To compare the false alarm rates of CUSUM to that of the proposed change detection scheme, the duration and the outcome of the FSS HT would need to be considered. The inclusion of false alarm rates for CUSUM in Figure 3 is there to justify the selected CUSUM threshold of $10^{2.157}$  for the comparison.

%We next consider the comparison between the proposed change detection scheme and the approach using a fixed sample size hypothesis test followed by a traditional change detection scheme. From Figure \ref{plotchangetp2}, it is clear that if the change time is known, a value of $\tau$ can be selected for the proposed change detection scheme which achieves a lower incorrect detection rate than the FSS HT. By cross referencing with Figure \ref{plotchangetp1}, it should be noted that for a given parameter value of $t>0$, over the range of change times over which it has a lower incorrect detection rate than the FSS HT, there is still an increased average detection delay caused by the initial state uncertainty cost. For example, for $\tau=$10^3$, the proposed change detector achieves lower incorrect detection rates for all change times of 35 through 65 when compared to the FSS HT; however, for a change time of 35, the proposed test has an average detection delay of 13.482, which is 51.4 \% larger than the average detection delay 
% However, the comparison between these two approaches is particularly meaningful for the case where the initial state is unknown and assumed to be small. For change times of $m=35$ and smaller, the proposed change detection scheme with $\tau=$0$ outperforms the FFS HT in terms of the incorrect detection rate.

\subsection{Analysis of Simulation Results} \label{subsec-analysis}

We next compare the results presented in Figures \ref{plotchangetp1} and \ref{plotchangetp2} with the performance bounds presented in the proof of Theorem \ref{theorem2}. As pointed out previously from Figure \ref{plotchangetp2}, once the initial state of the sequence is observed for long enough, each value of $\tau$ reaches a minimum probability of incorrect detection from initial state, and this minimum decreases as $\tau$ is increased. Eq. (\ref{upperboundallstates}) in Appendix~C shown below is an upper bound on the probability of incorrect detection from initial state. For $D=2$ with $f_0$ and $f_1$ defined  as in (\ref{PDF_test1}), $|\mathcal{S}_{-j}| = 1$. Since the PDFs are symmetric about $\mu/2$,
%(i.e. $f_0(x-\mu/2) = f_1(x+\mu/2)$), 
the performance of the test is the same regardless of whether $f_0$ or $f_1$ is the initial state. Without loss of generality we proceed by using (\ref{upperboundallstates}) and $j=0$ to calculate 
% the upper bound of the probability of incorrect detection, 
%From the previous arguments, for the PDFs $f_0$ and $f_1$ defined, the upper bound on the probability of incorrect detection from initial state is
\begin{align}
P(\mbox{incorrect detection from initial state}) \leq \max_{1 < m \leq n} \min_{s>0}  \frac{ \left( \mathbb{E}_1 \left[ \left( \frac{f_0(X_i)}{f_1(X_i)} \right) ^s \right] \right)^{m-1} }{\left( \frac{\tau}{C_1} - 1 \right)^{s} }. \label{pincanalysis}
\end{align}
Using the given PDFs, 
\begin{align}
\mathbb{E}_1 \left[ \left( \frac{f_0(X_i)}{f_1(X_i)} \right) ^s \right] = \exp \left( \frac{s(s-1)}{2} \frac{\mu^2}{\sigma^2}\right). \label{exps} 
\end{align}
We found the value of $C_1$ to be approximately 1.05 by taking an average of $\riskone{1}{n}{0}$ over $2 \leq n \leq 500$ for $\riskone{1}{500}{0}$ true over $10^3$ Monte Carlo trials. With this value of $C_1$ and (\ref{exps}), upper bounds (\ref{pincanalysis}) can be found for a given value of $\tau$ using two 1-dimensional searches. For values of $\tau$ of $10^1$, $10^3$, $10^5$, and $10^7$,  the probability of incorrect detection upper bounds were found to be \num{1.15e-1}, \num{1.00e-3}, \num{1.05e-5}, and \num{1.05e-7} respectively. Comparing values with Figure \ref{plotchangetp2}, (\ref{pincanalysis}) is a loose upper bound: incorrect detection probabilities for $\tau=10^1$ and $\tau=10^3$ are smaller than their respective upper bounds by factors of approximately \num{3e-3} and \num{2e-2} respectively.  This is expected, as the formulation of (\ref{pincanalysis}) considers only initial state uncertainty cost $\tau$, and ignores the exponential costs used to prevent incorrect detections.

In Theorem \ref{theorem2},  a minimum change time for initial state uncertainty risk to negligibly  increase delay was identified.  Eq. (\ref{prob_delay2}) is the probability that the initial state uncertainty risk is larger than expected risk $C_1$ given a change time  $m$, which incurs additional expected delay. It is desired to find the smallest $m$ such that the probability (\ref{prob_delay2}) is smaller than some small value $\alpha$. We compute 
\begin{align}
\sum_{i=1}^{m-1} \log \left( \frac{ f_0(X_i) }{ f_{\mathcal{S}_{-0}} (X_i) } \right) & = \sum_{i=1}^{m-1} \log \left( \frac{f_1(X_i)}{f_0(X_i)} \right), \nonumber \\
%& = \sum_{i=1}^{m-1} \log \left( \frac{\frac{1}{\sqrt{2 \pi \sigma^2}} \exp \left( \frac{-(X_i-\mu)^2}{2\sigma^2} \right)}{ \frac{1}{\sqrt{2 \pi \sigma^2}} \exp \left( \frac{-(X_i)^2}{2\sigma^2} \right)  } \right), \nonumber \\
%& = \sum_{i=1}^{m-1} \frac{X_i^2}{2 \sigma^2} - \frac{(X_i-\mu)^2}{2 \sigma^2}, \nonumber \\
%& = \sum_{i=1}^{m-1} \frac{2X_i \mu - \mu^2}{2 \sigma^2}, \nonumber \\
& = \frac{\mu }{\sigma^2} \sum_{i=1}^{m-1} X_i - \frac{(m-1)\mu^2}{2 \sigma^2}. \label{minctp1}
\end{align}
Noting that $X_{1,m-1}$ is distributed according to $\hypone{1}{m-1}{1}$, each of the $X_i$'s in (\ref{minctp1}) are IID $ \sim f_1$, 
%$\sum_{i=1}^{m-1} X_i$ is Gaussian with mean $(m-1)\mu$ and variance $(m-1)\sigma^2$, and 
(\ref{minctp1}) is Gaussian with mean $\frac{(m-1)\mu^2}{2\sigma^2}$ and variance $\frac{(m-1)\mu^2}{\sigma^2}$. Thus, (\ref{prob_delay2}) becomes
\begin{align}
1 - Q \left( \frac{ \log \left( \left( \frac{\tau}{C_m}-1 \right)(D-1)\right) - \frac{(m-1)\mu^2}{2\sigma^2} }{\frac{\sqrt{(m-1)}\mu}{\sigma}} \right) < \alpha. \label{probq}
\end{align}
%By selecting a small value of $\alpha$, (\ref{probq}) can be used to find the minimum change time $m$ for the initial state uncertainty risk to have a insignificant effect on the average detection delay. 
For $\alpha = 10^{-3}$, we obtain $C_m = 4$ by calculating the average values of $\riskk{m}{n}{0}{1}(\mbox{delay})$ over $2 \leq n \leq 500$ over $10^3$ Monte Carlo trials. For values of $\tau$ of $10^1$, $10^3$, $10^5$, and $10^7$, Eq. (\ref{probq}) yields minimum change times of 40.8, 59.2, 74.09, and 88.2, respectively. Interpolating linearly between data points of the average detection delay plot in Figure \ref{plotchangetp1}, the minimum change times for values of $\tau$ of $10^1$, $10^3$, $10^5$ and $10^7$ yield average detection delays of 9.261, 9.123, 9.115 and 9.106, respectively. The change detector for $\tau=0$ achieves an average detection delay of 8.90 over change times of 20 through 100, which is only slightly smaller, 
%Using Eq. (\ref{prob_delay2}) with $\alpha = 10^{-3}$ yields minimum change times that ensure average detection delay is 
within $4.02 \%$, $2.48 \%$, $2.38 \%$, and $2.28 \%$, respectively, of the above calculated minimum average detection delays.
% for $\tau$ of $10^1$, $10^3$, $10^5$, and $10^7$ respectively. 
Additionally, from Figure \ref{plotchangetp2}  it can be observed that, for each value of $\tau$, the incorrect detection rate achieved matches that calculated using (\ref{probq}), i.e., (\ref{prob_delay2})  identifies the duration of the initial transient period of the test.
%, during which increased average detection delay and an increased probability of incorrect detection are observed.

\subsection{Discussion} \label{subsec-summaryresults}

The above Monte Carlo simulations illustrate performance bounds and trade-offs presented in Theorem 2. Specifically, the initial state uncertainty cost $\tau$ can be chosen to trade off probability of incorrect detection from initial state according to a minimum change time and without appreciable average detection delay penalty. Furthermore, once this  minimum value of change time has been reached, initial state uncertainty risk decreases asymptotically to zero 
%for risks corresponding to the correct initial state 
and thus does not affect the trade-off between the probability of false alarm and the average detection delay.
%As such, once the parameters $a$, $b$, and $c$ are selected to satisfy a certain criteria regarding the trade-off between the average detection delay and the false alarm rate, the initial state uncertainty cost $\tau$ can be chosen without affecting this trade-off. 

%It is worth summarizing the selection of the parameters $a$, $b$, and $c$. 
In summary,  while joint optimization of parameter values $a$, $b$, $c$  and $\tau$ is complicated, we have instead shown that after a controllable initial transient period, the proposed change detection scheme can achieve detection delays close to that of CUSUM, which is optimal for the case where the initial state is known. Further, initial state uncertainty was considered  and it was found that that there is a range of change times where the proposed change detector outperforms FSS HT, which assumes full knowledge of the change time, in the form of the correct initial state detection probability. This can be attributed to the proposed change detector recursively tracking the minimum risk change time for each state pair $(r,s) \in \setsquare$, whereas  the FSS HT 
%assumes an IID sequence prior to the change occurring, and 
ignores the temporal behaviour of observed samples prior to the change. 

\section{Conclusions and Future Work}

A Bayesian change detection scheme has been proposed with exponential incorrect detection cost structure for a sequence of independent random variables that is known to start at one of $D$ equally likely possible PDFs and end at another.  Suitable parameter choices have been analytically shown to trade off the probability of incorrectly identifying the initial PDF against the test's false alarm probability and average detection delay and are illustrated by Monte Carlo simulation. Additionally, the simulations reveal that the proposed change detection scheme achieves an average detection delay close to that of the optimal change detector for the same ARL to false alarm and where the initial state is known. The analysis is restricted to the case where each of the PDFs have equal energy, which includes the case of different shifts in the mean. Generalizing the results to consider broader classes of PDFs is a subject of future investigation. 

\section{Appendix}
\subsection{Proof of Lemma~1:}
\begin{proof}
In (\ref{d}), both numerator and denominator take the form 
\begin{align}
\mathbb{E}_i [f_j(X)] & = \int_{- \infty}^{\infty} f_j(x)f_i(x)dx  \equiv \langle f_i, f_j \rangle \qquad \qquad i,j \in \mathcal{S}. \label{expectation}
\end{align}
%In (\ref{expectation}), $\langle \cdot , \cdot \rangle$ denotes, inner product. 
For $(j,k) \in \setsquare$,
\begin{align}
\dval{j}{j}{k} & = \frac{\mathbb{E}_j [ f_j(X)]}{\mathbb{E}_j[f_k(X)]}  = \frac{\langle f_j, f_j \rangle}{\langle f_j, f_k \rangle},   \nonumber \\
	& > \frac{\langle f_j, f_j \rangle}{ \sqrt{\langle f_j, f_j \rangle \langle f_k, f_k \rangle} }  = \frac{\sqrt{\langle f_j, f_j \rangle}}{\sqrt{\langle f_k, f_k \rangle}}
\end{align}
where the Cauchy-Schwarz inequality is used, and $\langle f_j, f_k \rangle > 0$ and $\langle f_j, f_j \rangle>0$ since both $f_j$ and $f_k$ are non-negative functions. The inequality is strict  since $f_j$ and $f_k$ are distinct and not linearly dependent. By letting $\langle f_j, f_j \rangle = \langle f_k, f_k \rangle \; \forall (j,k) \in \setsquare$, it follows straightforwardly that $\dval{j}{j}{k} > 1$. Similarly, $\dval{j}{k}{j} < 1$, noting that $\dval{j}{k}{j} = \left( \dval{j}{j}{k} \right)^{-1}$. 
\end{proof}

\subsection{Proof of Lemma~2:}
\begin{proof}
%
%It will be shown that if  $\mathbb{E} \big{[} \riskk{m}{n}{j}{k} P(X_{1,n}) \big{]}$  converges or diverges, then $\mathbb{E}[ \riskk{m}{n}{j}{k}]$ will also converge or diverge, respectively. 
Consider 
\begin{eqnarray}
\mathrm{Cov} \left[ \Big{(} \riskk{m}{n}{j}{k} P(X_{1,n}) \Big{)} , \frac{1}{P(X_{1,n})} \right] & = & \mathbb{E} \big{[} \riskk{m}{n}{j}{k} \big{]} - \mathbb{E} \big{[} \riskk{m}{n}{j}{k} P(X_{1,n}) \big{]} \mathbb{E} \left[ \frac{1}{P(X_{1,n})}) \right] \nonumber \label{cov1} \\
                                                                                                                                           & \leq & \sqrt{\mathrm{Var} \left[ \riskk{m}{n}{j}{k} P(X_{1,n}) \right] \mathrm{Var} \left[ \frac{1}{P(X_{1,n})} \right]}< \infty , \label{cauchyschwarz}
\end{eqnarray}
using the Cauchy-Schwarz inequality, 
and 
%is finite since the right side of (\ref{cov1})
%$\riskk{m}{n}{j}{k} P(X_{1,n})$ and $P(X_{1,n})^{-1}$ 
is finite by assumption that the PDFs have finite variance.
Applying Jensen's inequality 
%to convex function $1/x$ for $x>0$, 
\begin{align}
\mathbb{E} \left[ \frac{1}{P(X_{1,n})} \right] \geq \frac{1}{\mathbb{E} \big{[} P(X_{1,n})\big{]}}. \label{jensenxminusone}
\end{align}
 In (\ref{jensenxminusone}), $P(X_{1,n} = x_{1,n})$ is finite for any realization $x_{1,n}$ 
%of random sequence $X_{1,n}$, 
and thus the left and right hand sides of (\ref{jensenxminusone}) are a finite distance apart.
Applying (\ref{jensenxminusone}) to (\ref{cov1}) yields
\begin{align}
\mathbb{E} \big{[} \riskk{m}{n}{j}{k} \big{]} \geq \mathrm{Cov} \left[ \Big{(} \riskk{m}{n}{j}{k} P(X_{1,n}) \Big{)} , \frac{1}{P(X_{1,n})} \right] + \frac{\mathbb{E} \big{[} \riskk{m}{n}{j}{k} P(X_{1,n}) \big{]}}{\mathbb{E}[P(X_{1,n})]}. \label{cov2}
\end{align}
%The above inequality also follows from Jensen's inequality 
and
%in (\ref{jensenxminusone}), 
we  conclude
that $\mathbb{E}[ \riskk{m}{n}{j}{k} ]$ is a finite distance from $\mathbb{E} \big{[} \riskk{m}{n}{j}{k} P(X_{1,n}) \big{]}$.
% \left( \mathbb{E}[P(X_{1,n})] \right)^{-1}$. 

$\mathbb{E}_{\hypk{m}{n}{j}{k}}[P(X_{1,n})]$  can be calculated similarly to Eq.~(\ref{ERkn}) using (\ref{exphypkk})-(\ref{exphyponeone}) to be a product of $W_{m,n}(j,k)$ and a sum of geometric series terms, which along with Lemma \ref{lemma1}, can be shown to converge,
and thus $\mathbb{E}_{\hypk{m}{n}{j}{k}}[P(X_{1,n})]$ is a finite distance from $W_{m,n}(j,k)$. Taking expectation over $\hypk{m}{n}{j}{k}$, we conclude that $\mathbb{E}_{\hypk{m}{n}{j}{k}} \left[ \riskk{m}{n}{j}{k} \right]$ is finite. In a similar manner, using (\ref{cov2}) while taking expectation over $\hypk{m}{n}{p}{q}$ for $(p,q) \in \setsquare$ and $j \neq p$, $ \mathbb{E}_{\hypk{m}{n}{p}{q}} \left[ \riskk{m}{n}{j}{k} P(X_{1,n})\right]$ diverges, and thus $ \mathbb{E}_{\hypk{m}{n}{p}{q}} \left[ \riskk{m}{n}{j}{k}\right]$, a finite distance away, also diverges.
\end{proof}

\subsection{Proof of Theorem~2:}
\begin{proof}
\underline{Part (i):} 
From 
% arguments leading to 
Theorem~1, $\mathbb{E} \left[ \riskone{1}{n}{j} |  X_{1,n} \sim \hypone{1}{n}{j} \right] $  converges to a finite value. 
For incorrect detection at time $m$, there exists $\riskk{m}{n}{p}{q} < \riskone{1}{n}{j}$, for $j \in \mathcal{S}$, $(p,q) \in \setsquare$, and $p \neq j$ when $\hypone{1}{n}{j}$ is true.
%, and we expect $R_1^n \approx \mathbb{E} \left[ R_1^n | H_1 \mbox{ true}\right]$. 
 
Let $C_1 \stackrel{\Delta}{=} \mathbb{E} \left[ \riskone{1}{n}{j} |  X_{1,n} \sim \hypone{1}{n}{j} \right]  >0$.
% represent a risk threshold for $\riskk{m}{n}{p}{q} / {\cal Q}$ to be confident of incorrect detection. 
%Since the factor ${\cal Q}$ is common to all risks in the test, and has been factored out of all risks from here on. 
Since the parameter values $a$ and $c$ follow (\ref{cond1}) and $b<\infty$, $C_1 < \infty \; \forall n \geq 1$.
Using (\ref{Rkn}) under $\hypone{1}{n}{j}$, the all-$f_j$ hypothesis, an incorrect detection may only arise from  incorrect initial state $f_p$ when
\begin{align}
C_1 & > \tau \frac{ \sum_{\{ r \in \mathcal{S} | r \neq p \}} \prod_{i=1}^{m-1} f_r (X_i) }{\sum_{\{ r \in \mathcal{S}  \}} \prod_{i=1}^{m-1} f_r (X_i)} \mbox{, or equivalently,} \nonumber \\
%\frac{\tau}{C_1} & < \frac{ \prod_{i=1}^{m-1} f_p(X_i) + \sum_{\{ r \in \mathcal{S} | r \neq p \}} \prod_{i=1}^{m-1} f_r (X_i) }{\sum_{\{ r \in \mathcal{S} | r \neq p \}} \prod_{i=1}^{m-1} f_r (X_i)}, \nonumber \\
%\frac{\tau}{C_1} -1 & < \frac{\prod_{i=1}^{m-1} f_p(X_i)}{\sum_{\{ r \in \mathcal{S} | r \neq p \}} \prod_{i=1}^{m-1} f_r (X_i)}, \nonumber \\
%\log \left( \frac{\tau}{C_1}-1 \right) & < - \log \left( \sum_{\{ r \in \mathcal{S} | r \neq p \}} \frac{\prod_{i=1}^{m-1} f_r(X_i)}{\prod_{i=1}^{m-1} f_p(X_i)} \right), \nonumber \\
\log \left( \frac{\tau}{C_1}-1 \right) & < - \log \left( \prod_{i=1}^{m-1} \frac{f_j(X_i)}{f_p(X_i)} \right) - \log \left( 1 +\sum_{\{ r \in \mathcal{S} | r \neq p \mbox{ \small{and} } r \neq j \}} \prod_{i=1}^{m-1} \frac{f_r(X_i)}{f_j(X_i)} \right), \label{ineq_opp_bef}\\
%\log \left( \frac{\tau}{C_1}-1 \right) 
& < \sum_{i=1}^{m-1} \log \left(  \frac{f_p(X_i)}{f_j(X_i)} \right). \label{ineq_opp}
\end{align}
Eq.~(\ref{ineq_opp}) follows from (\ref{ineq_opp_bef}) since the omitted term (which only exists for $D>2$) is strictly less than zero.
%and thus whenever (\ref{ineq_opp_bef}) occurs, (\ref{ineq_opp}) must also have occurred. 
%This simplification does add some looseness to this condition; however, it can be shown that the expectation of the omitted term over $\hypone{1}{m-1}{j}$ goes to zero asymptotically as $m$ increases. 
%It is worth noting that the omitted term from (\ref{ineq_opp_bef}) only exists for $D>2$.
Let $Y_i \equiv  \log \left( \frac{f_p(X_i)}{f_j(X_i)} \right)$ and $Y \equiv \sum_{i=1}^{\overline{k}-1} Y_i$. Since the sequence
%the sequence of random variables 
$ \{ X_1, X_2, \ldots \}$ is IID, so is $ \{ Y_1, Y_2, \ldots \} $. For any $s>0$, the Chernoff bound yields
\begin{align}
P ( Y \geq \log ( {\tau}/C_1 - 1 ) ) &\leq  \mathbb{E}_j \left[ e^{sY} \right]  e^{-s \log \left( \frac{\tau}{C_1} - 1 \right) } \quad \label{chernoff}
\end{align}
where $\mathbb{E}_j \left[ e^{sY}\right] =  \prod_{i=1}^{\overline{k}-1} \mathbb{E}_j \left[ e^{sY_i}\right] $
%= \prod_{i=1}^{\overline{k}-1} M_{Y_i}(s)$ 
is the moment generating function of $Y$ and
%Since $Y_i$ for $1 \leq i < m$ is IID,  
%$M_Y(s)  = (M_{Y_i}(s))^{m-1}$.  
%Thus, 
%(\ref{chernoff}) can be rewritten as 
\begin{align}
P \left( Y \geq \log \left( \frac{\tau}{C_1} - 1 \right) \right) & \leq 
\frac{ \left( \mathbb{E}_j \left[ \left( \frac{f_p(X_i)}{f_j(X_i)} \right) ^s \right] \right)^{m-1} }{\left( \frac{\tau}{C_1} - 1 \right)^{s} }  ,  s > 0 \label{upperboundwithk}. \\
%\end{align}
%\begin{align}
%P \left( Y \geq \log \left( \frac{\tau}{C_1} - 1 \right) \right) 
& \leq  \max_{1 < m \leq n} \min_{s>0}  \frac{ \left( \mathbb{E}_j \left[ \left( \frac{f_p(X_i)}{f_j(X_i)} \right) ^s \right] \right)^{m-1} }{\left( \frac{\tau}{C_1} - 1 \right)^{s} } \label{upperboundmaxmin}
\end{align}
using recursively tracked change time, $m \in \{2, 3, \ldots , n\}$, and  
%according to hypothesis $\hypk{m}{n}{p}{q}$, can take any integer value from 2 to $n$. Considering all possible values of $m$, and 
minimizing over $s>0$.

The numerator of Eq. (\ref{upperboundmaxmin}) is positive for  $s > 0$, any pair $f_p$ and $f_j$, and any $m>1$.  Thus,
%on the probability of incorrect detection from incorrect initial state $f_p$ 
%Taking all possible incorrect initial states into account, 
%the probability of an incorrect detection occurring when $\hypone{1}{n}{j}$ is true is
\begin{align}
& P \left( \mbox{ incorrect detection }| X_{1,n} \sim \hypone{1}{n}{j} \right) \nonumber \\
= &  P \left( \bigcup_{\{ r \in \mathcal{S} | r \neq j \}} \mbox{ incorrect detection of initial state $f_r$} | X_{1,n} \sim \hypone{1}{n}{j} \right) \nonumber \\
\leq & \sum_{\{ r \in \mathcal{S} | r \neq j \}} \max_{1 < m \leq n} \min_{s>0}  \frac{ \left( \mathbb{E}_j \left[ \left( \frac{f_p(X_i)}{f_j(X_i)} \right) ^s \right] \right)^{m-1} }{\left( \frac{\tau}{C_1} - 1 \right)^{s} }. \label{upperboundallstates}
\end{align}
By choosing $\tau$ sufficiently large, Part (i) holds.

\underline{Part (ii):}

%The increase in delay caused by the initial state uncertainty risk will now be shown. With the exponential cost function chosen, 
%While increasing $\tau$ reduces the probability of incorrect detection, it also prevents detection of small change times. 
%We next show that 
%the correct-sided initial state uncertainty risk decreases exponentially towards zero as the recursively tracked change time increases and that other correct-sided risk components converge to finite non-zero values., so 
%the influence of $\tau$ on correct-sided test behaviour vanishes once the initial state uncertainty risk becomes sufficiently small. 
%The initial state uncertainty risk does not, however, change when the recursively tracked change time stays small. Consequently, for small change times, $\riskk{m}{n}{j}{k}$ may remain large enough to prevent it from becoming the smallest risk. This effect is analyzed as follows: 
%Similar to (\ref{recursive-Rnn}), the risk associated with recursively tracked change time $m$ can be computed (\cite{Falt17}, Eq. (3.19)). 
Denote the  risk of $\hypk{m}{n}{j}{k}$ excluding initial state uncertainty risk as
\begin{align} \label{gamma}
%\gamma \stackrel{\Delta}{=} 
\riskk{m}{n}{j}{k} (\mbox{delay}) \stackrel{\Delta}{=}  \riskk{m}{n}{j}{k} - \frac{\phi \tau}{(n-1)(D-1)+1}\frac{ \sum_{\{ r \in \mathcal{S} | r \neq j \} } \prod_{i=1}^{m-1} f_r(X_i)   }{  \sum_{\{ r \in \mathcal{S}  \} } \prod_{i=1}^{m-1} f_r(X_i)   }.
\end{align}
% where the expression on right-hand side differs from (\ref{wrong}) in that the change is from $f_0$ to $f_1$.

For the initial state uncertainty risk to vanish among all correct-sided risks, we require 
%the delay risk to exceed the average initial state uncertainly risk before a change occurs, i.e., 
%it must become smaller than $\gamma$. 
\begin{align}
C_m & > \tau \frac{ \sum_{\{ r \in \mathcal{S} | r \neq j \}} \prod_{i=1}^{m-1} f_r (X_i) }{\sum_{\{ r \in \mathcal{S}  \}} \prod_{i=1}^{m-1} f_r (X_i)},\label{partonecondition}
\end{align}
where $C_m \stackrel{\Delta}{=}  \mathbb{E} [ \riskk{m}{n}{j}{k} (\mbox{delay})  | \hypone{1}{n}{j} \mbox{~true}] / {\cal Q}$. Equivalently,
\begin{align}
\frac{\tau}{C_m} - 1 &< 
%\frac{ \prod_{i=1}^{m-1} f_j (X_i) }{ \sum_{\{ r \in \mathcal{S} | r \neq j \}} \prod_{i=1}^{m-1} f_r (X_i) }, \nonumber \\
\frac{1}{D-1} \left( \frac{ \prod_{i=1}^{m-1} f_j (X_i) }{ \frac{1}{D-1} \sum_{\{ r \in \mathcal{S} | r \neq j \}} \prod_{i=1}^{m-1} f_r (X_i) } \right), \nonumber \\
% &= \frac{1}{D-1} \left( \frac{f_j^{m-1}(X_{1,m-1} )}{\frac{1}{D-1}  \sum_{\{ r \in \mathcal{S} | r \neq j \}} f_r^{m-1} (X_{1,m-1}) } \right), \nonumber \\
 &= \frac{1}{D-1} \left( \frac{f_{j}^{m-1}(X_{1,m-1} )}{ f_{\mathcal{S}_{-j}}^{m-1} (X_{1,m-1}) } \right). \label{partoneconditioneqvexp}
 \end{align}
 Therefore,
\begin{align}
\log \left( \frac{\tau}{C_m}-1 \right) + \log \left( D-1 \right)&< \log \left( \frac{f_{j}^{m-1}(X_{1,m-1} )}{ f_{\mathcal{S}_{-j}}^{m-1} (X_{1,m-1}) } \right).
 \label{partoneconditioneqv}
\end{align}
Taking the expectation of (\ref{partoneconditioneqv}) conditioned on $\hypone{1}{n}{j}$  yields
\begin{align}
\log \left(\frac{\tau}{C_m} - 1  \right) &< \log \left( \frac{1}{D-1}\right) + \mathbb{E}_j \left[ \log \left( \frac{f_j^{m-1}(X_{1,m-1})}{f_{\mathcal{S}_{-j}}^{m-1}(X_{1,m-1})} \right) \right], \label{partoneconditionexp}
%k &> {1 + \frac{ \log \left( \frac{\tau}{  C_k } - 1 \right) }{ \mathbb{E}_0 \left[ \log \left( \frac{f_0(X_i)}{f_1(X_i)} \right) \right]}}. \label{logt}
\end{align}
where
\begin{align}
\mathbb{E}_j \left[ \log \left( \frac{f_j^{m-1}(X_{1,m-1})}{f_{\mathcal{S}_{-j}}^{m-1}(X_{1,m-1})} \right) \right] & =
\sum_{i=1}^{m-1} \mathbb{E}_j \left[ \log \left( \frac{ f_j(X_i)}{f_{\mathcal{S}_{-j}}(X_{i})} \right) \right] \nonumber \\
& = (m-1) \mathbb{E}_j \left[ \log \left( \frac{ f_j(X_i)}{f_{\mathcal{S}_{-j}}(X_{i})} \right) \right]  \geq 0. \label{KLdivergence}
\end{align}
In (\ref{KLdivergence}), $\mathbb{E}_j \left[ \log \left( \frac{ f_j(X_i)}{f_{\mathcal{S}_{-j}}(X_{i})} \right) \right] = D_{KL} \left( f_j || f_{\mathcal{S}_{-j}} \right)$. 
%which is greater or equal to zero with equality if and only if $f_j=f_{\mathcal{S}_{-j}}$. 
Using (\ref{KLdivergence}) and (\ref{partoneconditionexp}), (\ref{logt}) follows. 
%the average minimum change time for which the initial state uncertainty risk is smaller than other risk terms is
%\begin{align}
%m = \left \lceil{1 + \frac{ \log \left( \frac{\tau}{C_m} - 1 \right) + \log \left( D-1 \right)}{ D_{KL} \left( f_j || f_{\mathcal{S}_{-j}} \right) }} \right \rceil. \label{logt}
%\end{align}

%For the initial state uncertainty risk to increase the average detection delay by as little as possible, the probability of (\ref{partoneconditioneqv}) occurring should be as close to 1 as possible, i.e. it is desired to have

%To quantify the average detection delay, 
%The expected risks $\riskone{1}{n}{j}$ and $\riskk{m}{n}{j}{k}$ will be observed following a change at time $m$. 
Without loss of generality, we assume that the recursively tracked change time is the actual change time.
% i.e., the minimum risk change time for a change in the correct detection is the true change time. 
Detection delay is the difference $n-m$ after stopping to detect a change in the correct direction. Ignoring the influence of incorrect-sided risks, stopping occurs when 
\begin{align} \label{delay1}
\riskone{1}{n}{j} > \riskk{m}{n}{j}{k}.
\end{align}

%We proceed by taking the expectation of each side of (\ref{delay1}) over $X_{1,n} \sim \hypk{m}{n}{j}{k}$. 
The expectation of the right side of  (\ref{delay1}) is given by (\ref{ERkn}), and has been shown to converge to a finite value as $n$ gets large. Taking the expectation of the left side of  (\ref{delay1}),
\begin{align}
& \mathbb{E} [ \riskone{1}{n}{j} | X_{1,n} \sim \hypk{m}{n}{j}{k} ] \nonumber \\
= & \; \pi W_{m,n}(j,k) \Bigg{[} a^{n-m+1} \left( \sum_{i=2}^{m-1}  ( a \dval{j}{k}{j})^{m-i} \right)  + \sum_{i=m}^{n} a^{n-i+1} \dval{k}{j}{k}^{i-m} \nonumber \\
& + \sum_{\{ s \in \mathcal{S} | s \neq k \}} \Bigg{(} \sum_{i=2}^{m-1} a^{n-i+1} \dval{j}{s}{j}^{m-i} \dval{k}{s}{k}^{n-m+1}  + \sum_{i=m}^{n} a^{n-i+1} \dval{k}{j}{k}^{i-m} \dval{k}{s}{k}^{n-i+1} \Bigg{)} \nonumber \\
& + \; \mbox{other terms that do not contain parameter} \; a 
%& + \sum_{\{ r \in \mathcal{S} | r \neq j \}} \Bigg{(} c^n \dval{j}{r}{j}^{m-1} \dval{k}{r}{k}^{n-m+1} \nonumber \\
%& + \sum_{i=2}^{m-1} c^{i-1} \dval{j}{r}{j}^{i-1} \dval{k}{j}{k}^{n-m+1} + \sum_{i=m}^{n} c^{i-1} \dval{j}{r}{j}^{m-1} \dval{k}{r}{k}^{i-m} \dval{k}{j}{k}^{n-i+1} \nonumber \\
%& + \sum_{\{ s \in \mathcal{S} | s \neq j \mbox{ \small{and} } s \neq k \}} \Bigg{(}  \sum_{i=2}^{m-1} c^n  \dval{j}{r}{j}^{i-1} \dval{j}{s}{j}^{m-i} \dval{k}{s}{k}^{n-m+1} \nonumber \\ 
%& \hspace{100pt} +\sum_{i=m}^{n} c^n \dval{j}{r}{j}^{m-1} \dval{k}{r}{k}^{i-m} \dval{k}{s}{k}^{n-i+1} \Bigg{)} \nonumber \\
%& + c^{n-m+1} \sum_{i=2}^{m-1} c^{m-1} \dval{j}{r}{j}^{i-1} \dval{j}{k}{j}^{m-i} +  \sum_{i=m}^{n} c^n \dval{j}{r}{j}^{m-1} \dval{k}{r}{k}^{i-m} \Bigg{)} \Bigg{]}, 
\label{ER1n}
\end{align}
which contains terms that increase exponentially with base $a$ following the change at time $m$.

The average detection delay is the smallest $n-m>0$ satisfying the expectation of (\ref{delay1}) conditioned on $\hypk{m}{n}{j}{k}$. It can be shown that for $n>m$, $W_{m,n}(j,k) < \infty$, and observing (\ref{ERkn}), 
$\mathbb{E} [ \riskk{m}{n}{j}{k} | \hypk{m}{n}{j}{k}] < \infty$.
%following a change at time $m$, $\mathbb{E} [ \riskk{m}{n}{j}{k} | \hypk{m}{n}{j}{k} \mbox{ true}]$ converges to a finite value for $n>m$. 
For fixed $m$, the risk associated with initial state uncertainty is constant as $n$ increases. From (\ref{ER1n}),  $\mathbb{E} [ \riskone{1}{n}{j} | \hypk{m}{n}{j}{k}]$ has terms which increase as $a^n$. Since average detection delay is  the amount of time after $m$ for $\mathbb{E} [ \riskone{1}{n}{j} | \hypk{m}{n}{j}{k}]$ to surpass $\mathbb{E} [ \riskk{m}{n}{j}{k} | \hypk{m}{n}{j}{k}]$,  on average, the threshold $\tau$ increases detection delay by
\begin{align} \label{incavgdelay}
\frac{ \log \left( \tau \left(1 + \frac{1}{D-1} \exp ( (m-1) D_{KL}(f_j||f_{\mathcal{S}_{-j}}) ) \right)^{-1} \right) }{ \log \left( d_{min} \right)}, 
\end{align}
where $d_{\min}$ is determined by 
the largest value of $a$ satisfying 
(\ref{cond1}). Thus,  average detection delay increases  $O ( \log (\tau) )$, establishing (ii).
% Initial state uncertainty 

 \underline{Part (iii):} 
 Additionally, 
%as each of the $D$ distributions become more similar, as characterized by
under (\ref{dist_similar}), initial state uncertainty risk increases average detection delay, and it can be shown that $d_{min} \rightarrow 1$. Again, considering the increase in average detection delay from initial state uncertainty risk,  under (\ref{dist_similar}), the numerator of (\ref{incavgdelay}) approaches the constant $\log ( \tau (\frac{D-1}{D}) )$. Thus, delay from initial state uncertainty increases as $O( \log( d_{min} )^{-1})$.
\end{proof}

%It is well-known that the Kullback-Lieber divergence between any pdfs $f_j$ and $f_k$, $D_{KL} ( f_j || f_k ) = \int_{-\infty}^{\infty} f_j (x) \ln \frac{f_j (x)}{f_k (x)} dx \geq 0$ with equality iff $f_j = f_k$. Therefore, for $(j,k) \in \setsquare$, and thus $f_j \neq f_k$,
%\begin{eqnarray}
%1 < e^{D_{KL} ( f_j || f_k )} & \leq & \int_{-\infty}^{\infty} e^{\ln \frac{f_j (x)}{f_k (x)}} f_j (x) dx \nonumber \\
%& = & \int_{-\infty}^{\infty} \frac{f_j (x)}{f_k (x)} f_j (x) dx = \dval{j}{j}{k} ,
%\end{eqnarray}
%where we have used Jenson's inequality noting that $e^x$ is a convex function. Using a similar argument we obtain $d_0 > 1$.

%

%Once again considering Eq. (\ref{delay1}) for $f_1 \rightarrow f_0$, 

%the initial state uncertainty risk term will remain large since $d_0 \rightarrow 1$ when $f_1 \rightarrow f_0$. Equivalently , for changes in the opposite direction, $d_1 \rightarrow 1$ when $f_1 \rightarrow f_0$. As such, as $f_1 \rightarrow f_0$, the initial state uncertainty risk increases average detection delay as $O( log(d_0 d_1)^-1 )$.

%See  \cite{FaltARXIV2014}.

\bibliographystyle{IEEEtran}
\bibliography{IEEEabrv,refs}

% Generated by IEEEtran.bst, version: 1.14 (2015/08/26)
\begin{thebibliography}{10}
\providecommand{\url}[1]{#1}
\csname url@samestyle\endcsname
\providecommand{\newblock}{\relax}
\providecommand{\bibinfo}[2]{#2}
\providecommand{\BIBentrySTDinterwordspacing}{\spaceskip=0pt\relax}
\providecommand{\BIBentryALTinterwordstretchfactor}{4}
\providecommand{\BIBentryALTinterwordspacing}{\spaceskip=\fontdimen2\font plus
\BIBentryALTinterwordstretchfactor\fontdimen3\font minus
  \fontdimen4\font\relax}
\providecommand{\BIBforeignlanguage}[2]{{%
\expandafter\ifx\csname l@#1\endcsname\relax
\typeout{** WARNING: IEEEtran.bst: No hyphenation pattern has been}%
\typeout{** loaded for the language `#1'. Using the pattern for}%
\typeout{** the default language instead.}%
\else
\language=\csname l@#1\endcsname
\fi
#2}}
\providecommand{\BIBdecl}{\relax}
\BIBdecl

\bibitem{Falt2014}
J.~Falt and S.~Blostein, ``A {B}ayesian approach to two-sided quickest change
  detection,'' in \emph{IEEE ISIT}, 2014, pp. 736--740.

\bibitem{Falt2016}
------, ``Two-sided change detection under unknown initial state,'' in
  \emph{CISS}, 2016, pp. 418--423.

\bibitem{Blostein1991}
S.~Blostein, ``Quickest detection of a time-varying change in distribution,''
  \emph{IEEE Trans. Information Theory}, vol.~37, no.~4, pp. 1116--1122, 1991.

\bibitem{Moustakides1998}
G.~Moustakides, ``Quickest detection of abrupt changes for a class of random
  processes,'' \emph{IEEE Trans. Information Theory}, vol.~44, no.~5, pp.
  1965--1968, 1998.

\bibitem{Basseville93}
M.~Basseville and I.~Nikiforov, \emph{Detection of Abrupt Changes: Theory and
  Application}.\hskip 1em plus 0.5em minus 0.4em\relax Prentice-Hall, 1993.

\bibitem{Liang11}
J.~K. Liang and S.~Blostein, ``Performance evaluation of multiband multi-sensor
  spectrum sensing systems,'' in \emph{IEEE GLOBECOM}, 2011, pp. 1002--1007.

\bibitem{IEEE802.22Introduction}
C.~Stevenson, G.~Chouinard, Z.~Lei, W.~Hu, S.~Shellhammer, and W.~Caldwell,
  ``{IEEE 802.22: The first cognitive radio wireless regional area network
  standard},'' \emph{IEEE Communications Magazine}, vol.~47, no.~1, pp.
  130--138, Jan. 2009.

\bibitem{Page1954}
E.~S. Page, ``Continuous inspection schemes,'' \emph{Biometrika}, vol.~41, pp.
  100--115, June 1954.

\bibitem{Lorden1971}
G.~Lorden, ``Procedures for reacting to a change in distribution,'' \emph{Ann.
  Statist.}, vol.~42, no.~6, pp. 1897--1908, 1971.

\bibitem{Shiryaev1978}
A.~N. Shiryaev, \emph{Optimal Stopping Rules}.\hskip 1em plus 0.5em minus
  0.4em\relax Springer, 1978.

\bibitem{Moustakides1986}
G.~V. Moustakides, ``Optimal stopping times for detecting changes in
  distributions,'' \emph{Ann. Statist.}, vol.~14, no.~4, pp. 1379--1387, 1986.

\bibitem{Poor2008}
H.~V. Poor and O.~Hadjiliadis, \emph{Quickest Detection}.\hskip 1em plus 0.5em
  minus 0.4em\relax Cambridge University Press, 2008.

\bibitem{Zhao2010}
Q.~Zhao and J.~Ye, ``Quickest detection in multiple on-off processes,''
  \emph{IEEE Trans. Signal Proc.}, vol.~58, no.~12, pp. 5994--6006, Dec. 2010.

\bibitem{Shiryaev1963}
A.~Shiryaev, ``On optimum methods in quickest detection problems,''
  \emph{Theory Prob. App.}, vol.~13, no.~1, pp. 22--46, 1963.

\bibitem{Barnard1959}
G.~Barnard, ``Control charts and stochastic processes,'' \emph{J. Roy.
  Statistical Soc.}, vol.~21, pp. 239--271, 1959.

\bibitem{Dragalin1997}
V.~Dragalin, ``The design and analysis of 2-{CUSUM} procedure,''
  \emph{Communications in Statistics - Simul. Comput.}, vol.~26, no.~1, pp.
  67--81, 1997.

\bibitem{Barnett2001}
T.~Barnett, D.~Pierce, and R.~Schnur, ``Detection of anthropogenic climate
  change in the world's oceans,'' \emph{Science}, vol. 292, no. 5515, pp.
  270--274, 2001.

\bibitem{Hadjiliadis2005}
O.~Hadjiliadis, ``Optimality of the 2-{CUSUM} drift equalizer rules for
  detecting two-sided alternatives in the brownian motion model,'' \emph{J.
  Applied Prob.}, vol.~42, no.~4, pp. 1183--1193, 2005.

\bibitem{Hadjiliadis2006}
O.~Hadjiliadis and G.~Moustakides, ``Optimal and asymptotically optimal {CUSUM}
  rules for change point detection in the brownian motion model with multiple
  alternatives,'' \emph{Theory Prob. App.}, vol.~50, no.~1, pp. 131--144, 2006.

\bibitem{Hadjiliadis2008}
O.~Hadjiliadis and H.~Poor, ``On the best 2-{CUSUM} stopping rule for quickest
  detection of two-sided alternatives in a brownian motion model,''
  \emph{Theory Prob. App.}, vol.~53, no.~3, pp. 610--622, 2008.

\bibitem{Zarrin09}
S.~Zarrin and T.~J. Lim, ``Composite hypothesis testing for cooperative
  spectrum sensing in cognitive radio,'' in \emph{IEEE ICC}, 2009, pp. 1--5.

\bibitem{Pollak85}
M.~Pollak, ``Optimal detection of a change in distribution,'' \emph{Ann.
  Statist.}, vol.~13, no.~1, pp. 206--227, 1985.

\bibitem{Ritov90}
Y.~Ritov, ``Decision theoretic optimality of the {CUSUM} procedure,''
  \emph{Ann. Statist.}, vol.~18, no.~3, pp. 1464--1469, 1990.

\bibitem{Poor1998}
H.~Poor, ``Quickest detection with exponential penalty for delay,'' \emph{Ann.
  Statist.}, vol.~26, no.~1, pp. 2179--2205, 1998.

\bibitem{Poor1994}
H.~V. Poor, \emph{An Introduction to Signal Detection and Estimation}.\hskip
  1em plus 0.5em minus 0.4em\relax Springer-Verlag, 1994.

\bibitem{Falt17}
J.~Falt, ``Bayesian detection of a change in a random sequence with unknown
  initial and final distributions,'' Master's thesis, Queen's University, 2017.

\end{thebibliography}

\end{document}